\documentclass[11pt]{article}
\usepackage{placeins,graphicx,diagbox }

\usepackage{wrapfig} 
\usepackage{times} 
\usepackage{amsmath,tikz,pgfplots,mathtools,float}
\usetikzlibrary{arrows,shapes,snakes,shadows,positioning,automata,patterns}
\usetikzlibrary{trees,decorations.pathmorphing,decorations.markings}
\usepackage{amsmath}
\usepackage{subfigure}
\usepackage{multirow}
\usepackage{amssymb,latexsym}  
\usepackage{extarrows}
\usepackage{dsfont} 

\usepackage{amsmath,algorithm} 
\usepackage{amssymb}  
\usepackage{color}

\newtheorem{theorem}{Theorem}
\newtheorem{lemma}{Lemma}
\newtheorem{definition}{Definition}
\newtheorem{assumption}{Assumption}

\newtheorem{corollary}{Corollary}
\newtheorem{remark}{Remark}

\newcommand{\red}[1]{{\color{black}#1}}

\newcommand{\Real}{\mathbb{R}}

 \newcommand{\remove}[1]{}

\def\Real{\mathbb{R}}

\def\argmin{\mathop{\rm argmin}}

\def\us#1{{{\color{black}#1}}}
\def\uss#1{{{\color{black}#1}}}
\usepackage{verbatim}
\usepackage{ulem}
\author{Jinlong Lei \and Uday V. Shanbhag \thanks{ The authors  are with the
	Department of Industrial and Manufacturing Engineering, Pennsylvania
		State University, University Park, PA 16802  {\tt\small jxl800,udaybag@psu.edu}.
The work  has been partly supported  by NSF grant 1538605
	and 1246887 (CAREER). Additionally, a preliminary version of this work appeared in the Proceedings of the IEEE Conference on Decision and Control~\cite{Lei2017randomized}.}}

\date{}
\title{Asynchronous  Schemes for Stochastic and Misspecified Potential Games and Nonconvex Optimization}

\begin{document}
  \maketitle

\begin{abstract}
The distributed computation of equilibria and optima  has seen growing
{interest} in a broad collection of networked problems.
We consider the
computation of equilibria of convex stochastic Nash games characterized by a possibly
nonconvex potential function.  Since  any stationary point of the potential
function is a Nash equilibrium, there is
an equivalence between asynchronous  best-response (BR)  schemes applied on Nash game
and   block-coordinate descent (BCD) schemes implemented on the potential
function. We focus on two classes  of stochastic
Nash games:   (\textbf{P1}): A potential  stochastic  Nash
game, in which  each player solves a parameterized stochastic convex program; and
(\textbf{P2}): A misspecified    generalization,  where the player-specific
stochastic program is complicated by a parametric misspecification with  the
unknown parameter being the solution to  a stochastic convex  optimization.
In   both settings, exact proximal BR  solutions are generally unavailable in
finite time  since they   necessitate   solving   stochastic programs.
Consequently, we design two  asynchronous inexact  proximal BR  schemes to
solve problems (\textbf{P1}) and (\textbf{P2}), where in each iteration a
single  player  is randomly    chosen  to  compute an inexact  proximal BR
solution  (via stochastic approximation)  with delayed  rival   information
while    the other players {keep their strategies invariant}.  In the
misspecified regime  (\textbf{P2}),   each player possesses   an extra
estimate of the misspecified parameter  by  using a projected  stochastic
gradient (SG) algorithm with an increasing batch of sampled gradients.  By
imposing suitable conditions on the inexactness sequences, we prove that  the
iterates produced by both schemes   converge almost surely     to a
connected subset of the set of Nash equilibria.  {When the player  problems
are strongly convex, an inexact  pure  BR scheme is shown to be convergent.}
In effect, we provide what we believe is amongst the first randomized BCD
schemes for stochastic nonconvex  (but block-wise convex) optimization    with
almost sure  convergence properties.  We further show that the associated  gap function
converges to zero in mean.   These statements can be extended to allow for
accommodating weighted potential games and   generalized potential   games.
Finally,  we present preliminary numerics based on applying the proposed
schemes to  congestion control
and Nash-Cournot games.
 \end{abstract}

\section{Introduction}
  Nash games,  rooted in the  seminal work by \cite{nash50equilibrium}, have
seen wide applicability in  a broad range of engineered systems, such as power grids,
communication networks, transportation networks and sensor networks. In the $N-$player
Nash  game, each player maximizes a prescribed payoff
over a player-specific strategy set, given the rival strategies.    Nash's eponymous
solution concept,  Nash equilibrium (NE), requires that at an equilibrium, no
player can improve its payoff by unilaterally deviating from its equilibrium
strategy.  Potential games represent  an important   subclass of  Nash games, formally introduced in  \cite{monderer1996potential}, that
arise naturally in the modeling of many applications, ranging from
congestion control~\cite{kelly1998rate},  
routing   in communication networks \cite{altman2007evolutionary},  networked Cournot competition \cite{kannan2012distributed,yi2017distributed,yu2017distributed,jiang2017distributed},
 and a host of other control-theoretic problems \cite{ marden2009cooperative,scutari2006potential}.
We refer  the interested reader to the  survey by \cite{gonzalez2016survey} on potential games for  additional  references.    An interesting  extension  is the  class of {\em near potential games},  a class of games that are close to   potential games and  for which some learning dynamics (for instance, best-response, fictitious play,
 logit response) are discussed in \cite{candogan2013dynamics}. Additional decomposition algorithms are proposed in
\cite{facchinei2011decomposition} to solve {\em generalized} potential games  in which the  player-specific strategy set   depends on the strategies  selected by the other agents.

\noindent {\bf Motivation.}  While prior algorithmic efforts have considered
deterministic Nash games
(cf.~\cite{altman2007evolutionary,kannan2012distributed,yi2017distributed}),
there have {also been recent} attempts to contend with stochastic generalizations via
stochastic {gradient-response} schemes (cf.~\red{\cite{koshal2013regularized,  yousefian2013regularized, yousefian2016self, jiang2017distributed}}).
Yet, gradient-based schemes require {access to } {rival strategies}  after every update, often
undesirable in certain applications {such as} cellular networks. To this end, we
consider {\em inexact} {\bf proximal best-response} schemes applied to the
stochastic Nash game (or equivalently {\em inexact} {\bf block
coordinate-descent} schemes applied to the stochastic potential function) that
require  rival strategies   after   taking an inexact proximal BR step,
generally leading to  lower communication complexity (measured by the amount of rival information).  {Further,
in many regimes, the payoff functions are defined by parameters that are
unavailable, e.g., parameters of inverse-demand functions.}  {Accordingly}, we allow for
resolving {\bf parametric misspecification} in player  payoffs by equipping
each player with a simultaneous learning step.

\noindent {\bf Problems of interest.} We consider two classes of $N$-player  potential  stochastic Nash  games with players indexed by $i$ where $i\in \mathcal{N}\triangleq \{1,2,\cdots, N\}$.

 \textbf{(P1)}: {\em Potential  Stochastic Nash Games. }  Suppose the  $i$th player's strategy is denoted by $x_i$ with a strategy set
	  $X_i \subset    \mathbb{R}^{n_i}$,  implying that a feasible strategy $x_i$ satisfies  $x_i  \in X_i$,
and let $n \triangleq \sum_{i\in \mathcal{N}} n_i.$ Additionally, suppose   player $i$'s objective (or negative of the utility function)
	 is denoted by  $ f_i(x_i,x_{-i})$, which  depends  on its own strategy $x_i$ and on the {tuple} of
rival strategies  $x_{-i}\triangleq  \{x_j\}_{j \neq i}$.  {Suppose  $X$ and $X_{-i}$ are defined as $X   \triangleq \prod_{i=1}^N X_i$ { and }
$X_{-i}\triangleq   \prod_{j \neq i = 1}^N X_j,$  respectively. Given
rival strategies  $x_{-i}$, the  $i$th player is faced by the following     parameterized  stochastic  optimization problem:
 \begin{align} \label{problem1}\min_{x_i \in X_i}  \quad \ f_i(x_i,x_{-i})\triangleq  {\mathbb{E}_{\xi}}\left[
\psi_i(x_i,x_{-i};\xi(\omega)) \right],
\end{align}
where  {$\psi_i: X \times \mathbb{R}^d \to \mathbb{R}$} is a scalar-valued  function and the
random vector {is}  $\xi: { \Omega} \to \Real^d$  defined on the probability space $({ \Omega}, {\cal F}_x,
\mathbb{P}_x)$.  {Our interest lies in a subclass of Nash games,
	qualified as  \em potential},  characterized by a potential function    $P:X \to \mathbb{R}$
 such that for any $ i\in \mathcal{N}$ {and for any $x_{-i} \in X_{-i}$}:
\begin{equation}\label{def-potential0}
 P(x_i,x_{-i}) - P(x_i',x_{-i}) = f_i(x_i,x_{-i}) -f_i(x_i',x_{-i}),\quad  \forall x_i, x_i' \in X_i .
\end{equation}
Then  the Nash game, in which the $i$th player solves
\eqref{problem1} given $x_{-i}$, is called a  {\em potential   stochastic Nash game}.
We aim to compute  an NE, denoted by  $x^*=\{x_i^*\}_{i=1}^N$  such that  for any $ i\in \mathcal{N}$,  the following holds:
\begin{align} \label{minimal-optimal}
	f_i(x_i^*,x_{-i}^*) \leq f_i(x_i ,x_{-i}^*), \qquad \forall x_i  \in X_i.
\end{align}
 In other words,  a feasible strategy tuple  $x^* \in X$ is an Nash equilibrium  if no player
  can improve its payoff by unilaterally deviating from  its equilibrium  strategy $x_i^*$.

   \textbf{(P2)}: {\em Misspecified Potential Stochastic Nash  Games. }
 {The frequently used assumption in game-theoretic models is that each player has perfect knowledge  of the payoff function and is able to correctly forecast the choices of the other players.   However, as pointed out by \cite{kirman1975learning}     in the context of   Cournot oligopolies that  firms are, in general, imperfectly aware of their environment. Therefore they may have an imperfect knowledge of the payoff.  For example, players may employ  the misspecified  estimates  of  the  demand  function or the production capacity of their rivals.     \cite{kirman1975learning} introduced a learning  process where   firms update their conjectured  demand functions according to the observed data when the game is played repeatedly.   Subsequently, this notion was  formalized in a series of papers for resolving parametric misspecification  by}
 \cite{bischi2007oligopoly,bischi2008learning},
 \cite{szidarovszky2004global}, \cite{leonard1999nonlinear}, amongst others.
  Recent work has examined the development of coupled  {stochastic approximation (SA)} schemes for resolving misspecified
	stochastic optimization  \cite{jiang16solution,ahmadi2014resolution},  and stochastic Nash games \cite{jiang2017distributed}.
In this work, {we consider
  static stochastic Nash games  complicated by a parameteric misspecification $\theta^*$,}  in which the
$i$th player's  problem is  represented  as follows:
\begin{align} \label{problem1-mis}\min_{x_i \in X_i}  \quad \
f_i(x_i,x_{-i},\theta^*)\triangleq  {\mathbb{E}_{\xi}}\left[
\psi_i(x_i,x_{-i},\theta^*;\xi(\omega)) \right],
\end{align}
where $\theta^* \in \mathbb{R}^m$, $\xi: { \Omega} \to \Real^d$  is defined on the probability space   $({ \Omega}, {\cal F}_x,
\mathbb{P}_x)$,  and  $\psi_i: X \times \mathbb{R}^m\times  \mathbb{R}^d \to \mathbb{R}$  is a scalar-valued  function.
  For instance, in the context 	of Nash-Cournot games \cite{bischi2008learning}, $\theta^*$ may represent the slope and
	intercept of a linear inverse demand  (or price)  function; see the example of Nash-Cournot games with misspecified parameters in Section \ref{sec-numerics-formulation}. 	  We consider the estimation of   $\theta^*$   through solving  a  suitably defined convex stochastic  program:
\begin{align}\label{g_theta} \min_{\theta\in \Theta} g(\theta):= {\mathbb{E}_{\eta}}\left[g(\theta,\eta ) \right],
\end{align}
where $\Theta \in \mathbb{R}^m$ is a closed and convex set,   $\eta:\Lambda  \to \Real^p$  is defined on the probability space   $(\Lambda, {\cal F}_{\theta}, \mathbb{P}_{\theta})$,  and  $g: \Theta \times \mathbb{R}^p \to \mathbb{R}$  is a scalar-valued  function.  {
We still consider the class of  potential games and assume there exists  a  function
 $P(\cdot;\cdot ):X \times \Theta \to \mathbb{R}$  	 such that for any $ i\in \mathcal{N}$ and every $x_{-i} \in X_{-i}$:
\begin{equation}\label{def-potential}
 P(x_i,x_{-i};\theta^*) - P(x_i',x_{-i} ;\theta^*) = f_i(x_i,x_{-i};\theta^*) -f_i(x_i',x_{-i};\theta^*) ,  \quad  \forall x_i, x_i' \in X_i .
\end{equation} Then   we  refer to the  stochastic  Nash game    \eqref{problem1-mis} with  the misspecified parameter  $\theta^*$ being a solution to \eqref{g_theta}, as a misspecified potential stochastic Nash  game.}
	In such an instance, the 	problem of interest is to compute  the correctly specified Nash
	equilibrium, defined as follows, holds for $i \in \mathcal N$:
\begin{align} \label{minimal-optimal-mis}
	f_i(x_i^*,x_{-i}^*,\theta^*) \leq f_i(x_i ,x_{-i}^*,\theta^*), \qquad \forall x_i  \in X_i.
\end{align}
\noindent \red {{\bf Prior Research}. We now discuss some relevant prior research on stochastic Nash games, best-response schemes for Nash games, coordinate-descent schemes for  optimization problems,   and distributed {schemes for computing Nash equilibria}. }

\red{\em (i) Stochastic Nash games.}  Both \textbf{(P1)} and  \textbf{(P2)} belong to the class of static stochastic Nash games, as opposed to their dynamic variants, {representing a set} of Nash games in which the player objective functions
are expectation-valued. Tractable sufficiency conditions for existence of equilibria to such  games were   provided
in~\cite{ravat11characterization} in regimes where player objectives were
convex and could be nonsmooth.  {SA schemes for  a subclass of convex
stochastic Nash games were presented} under  Lipschitzian assumptions in
\cite{koshal2013regularized} via an iterative regularization technique  while
in~\red{\cite{yousefian2013regularized,yousefian2016self}}, Lipschitzian requirements were relaxed  by
utilizing a randomized smoothing approach. In both instances,  almost sure convergence
to an NE is guaranteed under suitable monotonicity requirements on the variational map.
Though  gradient-based schemes exist for solving stochastic Nash games  are
characterized by ease of implementation and lower complexity  in terms of each
player step, such schemes are characterized by the following properties: (i)  Players  {
require rival strategies after every gradient step, necessitating a  significant amount of   communication}; (ii)    Convergence theory is reliant on a relatively strong  monotonicity assumption on the gradient map;
  (iii) The schemes are synchronous. \uss{This motivates considering asynchronous best-response schemes, generally characterized by lower communication requirements and do not require a strong monotonicity property.}  \red{Finally, \cite{bernheim84rationalizable} and
 \cite{pearce84rationalizable} provide a characterization of
rational behavior in a non-cooperative game, where a rational player acts optimally given the decisions of her competitors; i.e. rational players will play a best-response strategy. \uss{This represents a cornerstone of much of the discussion in the context of such techniques and} \uss{further motivates the consideration of} such schemes. However, we emphasize that schemes based on other strategies (such as gradient-response,  etc.) are also of importance based on the setting being examined and the underlying assumptions imposed.}  \uss{Succinctly,} \red{we aim to} develop  implementable
asynchronous   proximal best-response  schemes, a class of techniques that can cope with delays and requires far less communication.

\red{\em (ii) Best-response (BR) and coordinate-descent schemes.}  In BR schemes,
each player  selects a BR strategy, given  current rival strategies
 \cite{fudenberg98theory,basar99dynamic}.   {\cite{nisan2011best} shows that for a class of games, it is best for
a player, given that the others are repeatedly   employing  best-response, to also repeatedly employ a  BR scheme.} There have been efforts to  extend BR  schemes to \uss{engineering applications}~\cite{scutari2010mimo}, where \uss{the} BR can be expressed in a
closed form.  Recently, \cite{pang2017two}  have proposed several variants  of the BR
schemes to solve the two-stage noncooperative games with risk-averse players.  Proximal BR schemes appear to have been  first discussed in \cite{FPang09},   where it is shown that the  set of fixed points of  the proximal BR map is equivalent to the set of Nash equilibria  when the player-specific problem is convex.   Additionally,
\cite{facchinei2011decomposition}  propose  several  regularized  Gauss-Seidel
 BR schemes for  generalized potential games  and show that  a  limit  point of the generated  sequence is an NE when each player's subproblem  is convex. Recall that  stationary
points of the potential function are    NE  of the original game
when the player-specific  problems are  convex.  Thus,  BCD methods may be employed for obtaining either stationary points of the potential function  (if nonconvex) or global minimizers (if convex), where  the coordinates are partitioned into several
blocks (each corresponding to a player in the associated Nash game) and at each
iteration,  a single block is  chosen to update while the other blocks remain unchanged.   Its original format
 dates back to \cite{d1959convex},   where blocks were updated {\em cyclically}. Convergence has been extensively
studied for both convex  and   nonconvex regimes  with either   differentiable or
nondifferentiable  objectives  \cite{tseng2001convergence}.  Notice
that the {\em asynchronous}  BCD schemes,  where at each epoch  a single block \us{is}
chosen,  are, in essence, identical  to asynchronous BR schemes.
\cite{nesterov2012efficiency}  considers  a {\em randomized}  BCD method that
performs a  gradient 	update on a randomly selected block.
  Extensions  to convex  nonsmooth  regimes have been studied extensively  \cite{richtarik2014iteration}
  while in nonconvex regimes,    
 ~\cite{xu2015block,davis2016asynchronous} {provide} convergence theory.
  \cite{ghadimi2016accelerated} proposed an accelerated gradient method to solve
nonconvex and possibly stochastic optimization,  while
\cite{ghadimi2016mini}  designed  a randomized stochastic
projected gradient  algorithm  to solve the constrained stochastic composite optimization.
  We summarize much of the prior work in  Table \ref{TAB-lit}, where we observe that there is no available a.s.
  convergence theory for (misspecified)  potential stochastic Nash games (or nonconvex stochastic programs) via BR (or BCD) schemes.

 \begin{table} [!htb]
 \centering
\newcommand{\tabincell}[2]{\begin{tabular}{@{}#1@{}}#2\end{tabular}}
\tiny
 \centering
     \begin{tabular}{|c|c|c|c|c|c|c|c|c|}
        \hline
 Problem& Literature     & Applicability  & stochastic &method  &rate &  a.s. & \tabincell{c}  {convergence\\ in mean} &misspecification
      \\ \hline
  \multirow{6}{*}{   Optimization}
 &
  \multirow{2}{*}{      \cite{xu2015block} } &    convex   &\multirow{2}{*}{  $\checkmark$} &  \multirow{2}{*}{ BSG}  &
$\mathcal{O}\left(1/k\right)$ &    \tabincell{c}  {--}  & -- &\multirow{2}{*}{ $\times$}
       \\  \cline{3-3}  \cline{6-8} & & nonconvex  & &  & -- & -- &  $\checkmark$   &
            \\ \cline{2-9}  &
\cite{ghadimi2016accelerated}  & nonconvex   & $\checkmark$  &\us{gradient-based}&$\mathcal{O}(1/\sqrt{k})$&-- &$\checkmark$ & $\times$
    \\ \cline{2-9}  &
 \cite*{ghadimi2016mini} & nonconvex   & $\checkmark$ &\us{gradient-based} &$\checkmark$ &-- &$\checkmark$ & $\times$
 \\ \cline{2-9}   & \cite{davis2016asynchronous}  & nonconvex &$\times  $ &\tabincell{c}   {gradient-based\\ (asynchronous)}
 &  -- &\tabincell{c}{cluster point \\ is an NE} & --  & $\times$
      \\   \cline{2-9}   & {\bf This work}  &\tabincell{c}  { nonconvex \\ block-convex} &$\checkmark$&\tabincell{c}  {RBCD-based\\ (asynchronous)}  &-- & $\checkmark$ &$\checkmark$&
 \tabincell{c} {      Section  \ref{sec:randomized}: $\times$  \\  Section   \ref{sec:misp}: $\checkmark$   }\\
      \hline \hline
      \multirow{9}{*}{    \tabincell{c}  {Non-potential \\Game}}   &  \cite*{jiang2017distributed} & monotone   &$\checkmark$ & gradient-based  &$\mathcal{O}\left(1/k\right)$ & $\checkmark$ &--& $\checkmark$
      \\ \cline{2-9}
      &  \cite*{koshal2013regularized} & \tabincell{c}  { monotone, \\Lipschitz}   &$\checkmark$ &gradient-based&--& $\checkmark $ &--& $\times$
  \\ \cline{2-9}
      &  \cite*{yousefian2013regularized} & \tabincell{c}  { monotone, \\non-Lipschitz}   &$\checkmark$ &gradient-based&--& $\checkmark$ &--& $\times$
       \\ \cline{2-9} &       \red{ \cite*{yousefian2016self} } &\tabincell{c}{\red{strongly monotone,}\\ \red{non-Lipschitz}} &\red{$\checkmark$}&\red{gradient-based }&\red{$\mathcal{O}(1/k)$}& \red{$\checkmark$}  &\red{$\checkmark$} &\red{$\times$}
\\ \cline{2-9}
       &             \cite*{pang2017two} & \tabincell{c}  { contractive \\BR maps}
        & $\checkmark$ & {BR-based}  & linear &$\checkmark$   & -- & $\times$
       \\ \cline{2-9}
       &             \cite*{lei2017synchronous} & \tabincell{c}  { contractive \\BR maps}
        & $\checkmark$ & BR-based   &linear  &$\checkmark$&$\checkmark$ & $\times$
            \\ \hline \hline
           \multirow{3}{*}{    \tabincell{c}  {Potential \\Game}}   &   \cite*{facchinei2011decomposition} & \tabincell{c}  {generalized  } &
               $\times$ & \tabincell{c}  {BR-based \\(Gauss-Seidel)} & --  &\tabincell{c}{cluster point \\ is an NE}
               &--& $\times$
      \\   \cline{2-9}   & {\bf This work}  &\tabincell{c}  { player-specific \\ convex} &$\checkmark$&\tabincell{c}  {BR-based\\ (asynchronous)}  &-- & $\checkmark$ &$\checkmark$&
 \tabincell{c} {      Section  \ref{sec:randomized}: $\times$  \\  Section   \ref{sec:misp}: $\checkmark$   }
       \\ \hline
      \end{tabular}
      \vskip 3mm
        \centering \caption{A list of some recent research papers on  Nash games   and (non)convex optimization. }\label{TAB-lit}
 {\scriptsize   \red{The   3rd column ``Applicability"  specifies the     major conditions   required by  the  problem   studied in the literature.}  In column 4 ``stochastic" (or column 9 ``misspecification"), $\checkmark$ means that  the studied  problem in the literature  is stochastic (or with parametric  misspecification),  while $\times$ implies that the studied  problem   is  deterministic (or without parametric misspecification); In columns 6-8,  $\checkmark$implies that the reference has studied the corresponding convergence property, while dash  -- implies that  the literature has not studied this   property. }
\end{table}

\red{{\em (iii)  Distributed  computation of  Nash equilibria.}
Recently, there has been  an interest  in considering settings where players  compute their strategies in a distributed sense. 
\cite{grammatico2015decentralized} and \cite{parise2015network}
considered the networked aggregative games (where player payoffs are coupled via the aggregate strategy) with quadratic payoffs and proposed decentralized schemes for NE computation. In this setting, the communication graph is, in essence, the payoff dependence network. However, as has been pointed out in \cite{marden2007learning}, a player is inherently limited in that it can communicate with (and be observed by) a few players in large-scale networked systems.  Thus,
in recent work, it is assumed that players form a (connected)
communication network (which is not necessarily the same as the players' payoff
dependence network) and exchange messages with the neighboring players to
obtain estimates of unobserved information.  In such a setting,  \cite{koshal2016distributed}   developed consensus-based distributed algorithms for aggregative games while \cite{salehisadaghiani2016distributed} presented an asynchronous gossip-based algorithm. In a similar vein,  \cite{zhu2016distributed}  designed similar distributed
algorithms for a class of generalized convex games while  \cite{tatarenko2018geometric} considered similar settings under a strong monotonicity requirement on the concatenated gradient map. \cite{shi2017lana} designed a linearized
ADMM-like scheme  while continuous-time schemes were presented by \cite{ye2017distributed} in which unobservable decisions are learnt via a consensus-based approach. Finally,  \cite{lei2018distributed} have developed distributed counterparts of (inexact) best-response and gradient-response for a range of stochastic Nash games.}

  {\bf Contributions:}   		
 In \cite{shanbhag2016inexact} and \cite{lei2017synchronous},   rate statements and iteration complexity bounds are provided for inexact  proximal BR schemes  for stochastic Nash
games under a contractive requirement on the proximal  BR map. However,  even asymptotic guarantees are unavailable without such an assumption.   Motivated by this gap, we  aim to design  convergent   implementable
asynchronous   BR  schemes such that at each
epoch,  a single player updates its strategy  while  the other  players keep their strategies invariant.
  In our settings,  each  player-specific subproblem involves solving a stochastic program   whose
 exact solution is generally unavailable  in finite time,  necessitating inexact solutions.
 Accordingly, we  propose two classes of  asynchronous    inexact  proximal BR  schemes  to compute NE of
  problems \textbf{(P1)} and \textbf{(P2)},  and  make the following contributions: \noindent  {\bf (i).}   In Section II, 	we propose an asynchronous    inexact  proximal  BR  scheme  to solve \textbf{(P1)}.
     In each iteration, a single agent is randomly chosen to inexactly solve  a stochastic optimization problem,  given the delay-afflicted     rival strategies  via an SA scheme}.
 By imposing  suitable conditions on   \textbf{(P1)} and on the inexactness sequence,  in a regime that allows for  uniformly bounded delays,  we prove that  the iterates  converge a.s. to a connected subset of the
 set of   Nash equilibria and that  the  gap function  converges   in mean to zero.
Extensions are provided to  the generalized stochastic
potential games (with coupled strategy sets) and  the weighted potential games.   {We further prove that asynchronous    inexact  pure  BR schemes are convergent if player-specific problems are strongly convex.}
 \noindent {\bf (ii).} In Section III, we   extend the regime to contend with
the  misspecified stochastic Nash game \textbf{(P2)} where every player updates its equilibrium strategy and its belief regarding the misspecified parameter (via variable sample-size SA schemes),  given   rival strategies afflicted  by delays.
Asymptotic guarantees analogous to Section II are provided and we additionally
show that the belief regarding the misspecified
parameter  converges  a.s. to its true counterpart.
 \noindent (iii)   We provide some preliminary numerics on congestion games and Nash-Cournot games in   Section IV,  and
conclude  the paper in Section V.

  {\bf Notations:} When referring to a vector $x$, it is assumed to
be a column vector while $x^T$ denotes its transpose. Generally, $\|x\|$ {denotes}
the Euclidean vector norm, i.e., $\|x\|=\sqrt{x^Tx}$.   For a  nonempty closed convex   set $X \subset \mathbb{R}^m$,
we use $\Pi_X[x]$ to denote the Euclidean projection of
a vector $x\in \mathbb{R}^m$ on   $X$, i.e., $\Pi_X[x]=\min_{y \in X}\|x-y\|$.
  We write \textit{a.s.} as the abbreviation for ``almost surely''. We use $\mathbb{E}[{z}]$ to denote the
unconditional  expectation of a random variable~$z$.
For a real number $x $,   the floor function  $\lfloor  x \rfloor $ denotes the  largest integer smaller  than $x$.
   We use   $[A]_{i,j}$  to denote  the $(i,j)$-th entry of the matrix $A.$

   \section{Asynchronous  Inexact Proximal Best-Response  Schemes} \label{sec:randomized}
  In Section~\ref{Sec:IIA}, we  propose an asynchronous  inexact proximal  BR  scheme  to  compute an equilibrium of  the
   stochastic potential  game \textbf{(P1)}.  Then  in  Section~\ref{Sec:IIB} we
	   introduce    some basic assumptions, based on which, we proceed to prove
the almost sure  convergence  and  convergence in mean of the  generated sequence
to a Nash equilibrium in Section~\ref{Sec:IIC}. In Section~\ref{Sec:IID}, we  discuss  some possible extensions
	including the generalized  potential  games    allowing  for coupled
		strategy sets,  and weighted potential games.  {Finally, in Section \ref{sec:Br}, we show that in a delay-free regime,  the   asynchronous  inexact   pure BR  scheme  (i.e. without the proximal term) is  convergent when the player-specific problem is strongly convex.}

   \subsection{Algorithm Design}\label{Sec:IIA}
  In standard potential games, a natural approach is  an   asynchronous BR    method  where in each  iteration, one player updates its  strategy by solving problem \eqref{problem1},  given its   rival strategies,  referred to as the {\em best-response}	  problem.  {However, best-response schemes do not always lead to
convergence to Nash equilibria. In fact, even in potential
games where the potential function is player-wise convex, such
convergence does not follow; see \cite{facchinei2011decomposition} for a
simple  counterexample. Accordingly,
 \cite{facchinei2011decomposition}  propose a regularized  BR scheme in which  each player's objective is modified
   by adding a quadratic proximal term,   and prove  its convergence. Our research has been motivated by considering stochastic
generalizations which do not follow immediately.   Yet another
reason for using proximal  BR schemes may be  drawn from   the  ``momentum
behavior"  in economics, e.g. in the  investment problems \cite{grinblatt1995momentum}, the players
may want to optimize their objective  while staying  close to their
previous values.}
We then  define the player $i$'s  proximal best-response problem as follows for some $\mu_i>0$:
\begin{equation}\label{proximal}
 T_i(x) \triangleq \argmin_{y_i \in X_i}  \left[ \mathbb{E}\left[\psi_i(y_i, x_{-i};\xi(\omega)) \right]  +\frac{\mu_i}{2} \|y_i-x_i\|^2\right].
\end{equation}
Since $T_i(x)$, the minimizer of a  stochastic problem \eqref{proximal},   is generally unavailable in finite time, we
utilize Monte-Carlo sampling schemes in  obtaining inexact solutions ~\cite{shapiro09lectures}.

 {We assume that each player $i$ always knows its current strategy, while is not immediately  aware of  rival
strategies. Instead,   its knowledge of each rival  strategies may be afflicted by a
rival-specific random delay,  (see  \cite{fudenberg2014delayed}).}
We   now propose   an asynchronous   inexact   proximal best-response    scheme   (Algrithm \ref{algo-randomized}) to compute  an NE of   this game.   {At time $k$, player $i$'s  strategy  $x_i(k) \in\mathbb{R}^{n_i}$   is an estimate for its equilibrium strategy $x_i^* $  and player $i$ has access to   delayed  rival strategies $y^i(k) \triangleq (x_1(k-d_{i1}(k)),\cdots,x_N(k-d_{iN}(k))  )$,  where  $d_{ij}(k)$ denotes the  delay  associated with  player $j$'s
information, and $d_{ii}(k)=0$. }  The scheme  is  defined as follows.   At   iteration $k\geq 0,$ randomly pick  a single $ i$ from $ \mathcal{N} $ with probability   $\mathbb{P}(i_k=i)=p_i>0$.
 If $i_k=i,$  then player $i$ is chosen to  initiate  an update  by computing  an inexact proximal BR solution to problem \eqref{proximal} characterized by \eqref{randomized-alg1}.
 {We   impose conditions on the inexactness sequence    $ \{ \varepsilon_i(k)\}_{k\geq 1} $
when we proceed to   investigate  the convergence properties.}

 \begin{algorithm} [H]
\caption{{Asynchronous inexact proximal best-response scheme }}\label{algo-randomized}
 Let $k:=0$, $ x_{i,0} \in X_i$  for $ i\in \mathcal{N}$.
		 Additionally $0 < p_i < 1$ for $ i\in \mathcal{N}$ such that
$\sum_{i=1}^N  p_i = 1.$
\begin{enumerate}
\item[(S.1)] Pick $i_k=i \in \mathcal{N} $ with probability $p_i$.
\item[(S.2)] If $i_k=i$,  then    player $i$  updates     $ x_i(k+1) \in X_i$  as follows:
\begin{equation}\label{randomized-alg1}
x_i(k+1) :=  T_i(y^i(k))+ \varepsilon_i(k+1),
\end{equation}
 {where  $  \varepsilon_i(k+1) $ denotes the inexactness employed by player $i$ at time $k+1$.}
Otherwise,     $ x_j(k+1):= x_j(k)$ if $j  \notin  i_k$.
 \item[(S.3)]  {If $k > K$, stop; Else, $k:=k+1$} and return to (S.1).
\end{enumerate}
\end{algorithm}

 {	

\begin{remark}
 In fact,  in practical game-theoretic problems,
players take actions  in an asynchronous manner since there might
not exist a global coordinator  to ensure that  players update simultaneously.    The condition that
$\mathbb{P}(i_k=i)=p_i>0$ with $\sum_{i=1}^N p_i=1$ accommodates
the  Poisson model employed by \cite{aysal2009broadcast} and \cite{boyd2006randomized}  as   a special case.   For $i \in \cal{N}$, player $i$  is activated  according to   a local Poisson clock,  which  ticks according to a Poisson process with rate $\varrho_i>0$.
Suppose that there is a virtual global clock which ticks whenever
any of the local Poisson clocks tick.  Assume that the local
Poisson clocks are independent, then  the global clock ticks
according to a Poisson process with rate $\sum_{i=1}^N\varrho_i$.
Let  $Z_k$ denote the time of the $k$-th tick of the global
clock.  Since the  local Poisson clocks are independent,  with
probability one,   there is  a single player  whose Poisson clock  ticks at time $Z_k
$  with probability  $\mathbb{P}(i_k =i)={\varrho_i \over
\sum_{i=1}^N\varrho_i} \triangleq p_i.$ Further, the memoryless
property of the Poisson process indicates  that $\{i_k \}_{k\geq 0}$ is an independent and identically distributed (i.i.d.) sequence.
\end{remark}

}


\subsection{Assumptions  and Preliminary Results}\label{Sec:IIB}

 For notational simplicity, let $\xi$ denote $\xi(\omega)$ throughout the
 paper. We begin by imposing assumptions on $X_i$, $f_i$, $\psi_i$, and
 on the second moments of $\psi_i$.

\begin{assumption}\label{assump-play-prob} ~ Let the following hold.\\
(a) For every $i\in \mathcal{N}$,  the feasible set $X_i $ is  closed,  compact, and    convex;\\
(b)  For every $i\in \mathcal{N},$  $f_i(x_i,x_{-i})$ is convex and   continuously  differentiable  in $x_i \in X_i$ for every
$x_{-i} \in X_{-i}$.  {In addition,  there exists a Lipschitz constant  $L_i >0$ such that   the following holds:
$$\| \nabla_{x_i} f_i(x )-\nabla_{x_i} f_i(x')\| \leq L_i \| x-x'\|\quad \forall x,x'\in X  ;$$}
(c) For every $i\in \mathcal{N},$  all $x_{-i} \in X_{-i}  $, and any  $\omega \in  \Omega$, $  \psi_i(x_i, x_{-i};\xi(\omega))$ is differentiable  in
$x_i$ over {an open set containing} $X_i$ such that  $ \nabla_{x_i} f_i(x_i,x_{-i})=  \mathbb{E}_{\xi}[\nabla_{x_i} \psi_i(x_i, x_{-i};\xi )]$;\\
(d) For   every $ i\in \mathcal{N}$  and any  $x \in X $, there exists a constant $M>0$ such that
 $\mathbb{E}_{\xi}[\|\nabla_{x_i} \psi_i(x_i,x_{-i};\xi)\|^2] \leq  M^2.$
\end{assumption}
It is seen that (c) and (d) pertain to the existence of a
conditionally unbiased stochastic oracle and the boundedness
of the conditional second moment of the sampled gradient
generated by this oracle.   Next, we assume the existence of a continuously differentiable potential
function for the Nash game of interest.
\begin{assumption}[{\bf Potential function}]\label{assp-potential} ~There exists a  potential  function
$P: X\to \mathbb{R}$ that is continuously differentiable over an open set containing $X$    such that
 {for any $ i\in \mathcal{N}$   and  any $x_{-i} \in X_{-i}$, equation
 \eqref{def-potential0}  holds.} 	 \end{assumption}

  	 We next make some assumptions on the  delays as well as on
the inexactness sequences   $\{ \varepsilon_i(k)\}$ utilized  in Algorithm   \ref{algo-randomized}.
 {We denote  the $\sigma$-field of the entire information  used by  Algorithm \ref{algo-randomized} up to  (and including)  the update of  $x(k)$ by $\mathcal{F}'_k$, and the    $\sigma$-field generated from     $\mathcal{F}'_k$  and the   delays  at  time $k$  by   $\mathcal{F}_k\triangleq \sigma \big\{ \mathcal{F}_k',    d_{ij}(k), i,j\in \mathcal{ N}  \big\}$. We will  formally  define  $\mathcal{F}'_k$ after introducing  the SA scheme \eqref{sa-inner}.}
 \begin{assumption}[{\bf Delay  and inexactness sequences}] \label{assp-noise}
The following hold: \\
(a)   $i_k$ is independent of $\mathcal{F}_k$   for all $k\geq 1;$ \\
(b) for any $i \in \mathcal{N},$ the noise term $\{ \varepsilon_i(k)\}$
 satisfies the following condition:
 $$ \sum_{k=1}^{\infty}  \mathbb{E} \left[    \| \varepsilon_i(k+1)\|^2        \big | \mathcal{F}_k\right]  < \infty,~a.s.,
~ \textrm{and}~\sum_{k=1}^{\infty}  \mathbb{E} \left[    \| \varepsilon_i(k+1)\|        \big | \mathcal{F}_k\right]  < \infty,~a.s.;   $$
   (c)   there exists a positive integer $\tau$
such that for any $i,j\in \mathcal{N}$ and any $k\geq0,$ $d_{ij}(k) \in \{0,\cdots, \tau\}$.
\end{assumption}

By  Assumption \ref{assump-play-prob},  it is clear   that  $T_i(x) $ defined by
\eqref{proximal} requires solving a strongly convex stochastic program.
Thus, an  approximation of the solution  to  the problem  \eqref{proximal} {with $x=y^i(k)$}, characterized by
\eqref{randomized-alg1}, can be computed via the standard SA algorithm defined  as follows   for   $t = 1, \hdots, j_i(k)$:
\begin{equation} 
\begin{split}
& z_{i,t+1}(k)   := \Pi_{X_i} \Big[   z_{i ,t}(k)  -   \gamma_{i,t}   \left[\nabla_{x_i}  \psi_i(  z_{i,t} (k) ,y^i_{-i}(k);\xi_{i,t}(k)) + \mu_i (   z_{i ,t}(k)  -x_i(k))\right] \Big], \label{sa-inner}
 \end{split}
\end{equation}
 where  $ \gamma_{i,t}=\frac{1}{\mu_i (t+1)}$ and  { $ z_{i ,t}(k) $   denotes  the  estimate  of the proximal BR  solution  $T_i(y^i(k))$ at $t$-th inner step of the SA scheme  \eqref{sa-inner} with the intinal value $z_{i,1}(k) =x_i(k)$. Set  $x_i(k+1)   = z_{i, j_i(k) }(k).$}   Define $\xi_i(k)\triangleq ( \xi_{i,1}(k),\cdots, \xi_{i,j_i(k)}(k)),$ and   {$\mathcal{F}_k'\triangleq\sigma\{ x(0),i_l,\xi_{i_l}(l), d_{i_l,j}(l), 0 \leq l
  \leq k-1, j\in \mathcal{N}\}$.}
  Then by Algorithm \ref{algo-randomized},  $x(k)$  is adapted to $\mathcal{F}_k'$
  and   $y^i(k)$  is adapted to $\mathcal{F}_k$.  As a result,   $T_i (y^i(k) )$ is   adapted to $\mathcal{F}_k$  by   \eqref{proximal}.  Analogous  to   Lemma 3 in  \cite{lei2017synchronous},     the following result   holds for  the  SA scheme \eqref{sa-inner}.

\begin{lemma}\label{potential-lem-rate-sa} ~Let Assumption~\ref{assump-play-prob} hold.
Consider the asynchronous inexact proximal best-response scheme given by   Algorithm \ref{algo-randomized}.
Assume  that the random variables $\{ \xi_{i,t}(k)\}_{1 \leq t \leq j_i(k)}$ are i.i.d., and   that  for any $  i \in \mathcal{N}$  the random vector    $ \xi_i(k)$ is independent  of $\mathcal{F}_k$.   Then we have the  following for any  $t: 1\leq t \leq j_i(k)$.
 $$ \mathbb{E}[\|  z_{i,t}(k)-  T_i(y^i_k)\|^2\big | \mathcal{F}_k] \leq   {Q}_i/{(t+1)} \quad a.s.,$$
where $Q_i \triangleq  \frac{2M^2}{ \mu_i^2}+2D_{X_i}^2$   and    $ 	D_{X_i} \triangleq \sup\{ d(x_i,x_i'): x_i,x_i' \in X_i \}.$ \end{lemma}

 \begin{remark}\label{rem1}~
\noindent(i) Let  {${\bf 1}_{[i_k=i]}$} denote the indicator function
of the event  $i_k = i$, defined as  ${\bf 1}_{[i_k=i]}=1  ~{\rm if ~~}   i_k=i,$ and $=0,$ otherwise.
Define $  \Gamma_{i,0} \triangleq 1$ and
$\Gamma_i(k) \triangleq 1+\sum\limits_{t=0}^{k-1} {\bf 1}_{\{ i_t=i\}}$
for all $k \geq 1$.  Then the computation of $\Gamma_i(k)$ merely uses player $i$'s local information  and   for every $\omega \in \Omega $,
there exists  a sufficiently large $\tilde{k}(\omega)$  that is possibly contingent on the sample path $\omega $ such that  for  any $i\in \mathcal{N}: $
\begin{align}\label{lem-step-size}
\Gamma_i(k) \geq     { k p_i \over 2}+1 \quad \forall k \geq \tilde{k}(\omega).
\end{align}
 The proof can be found in Lemma 7 of \cite{koshal2016distributed}.

\noindent (ii)  Set  $j_i(k)\triangleq \left \lfloor  \Gamma_i(k)^{2(1+\delta)}\right \rfloor$   {for some positive $\delta>0$}
and  $x_i(k+1)   = z_{i, j_i(k) }(k).$
Then by   Lemma \ref{potential-lem-rate-sa},   we have that $ \mathbb{E}[\|x_i(k+1)-  T_i(x(k))\|^2\big | \mathcal{F}_k] \leq   \frac{Q_i}{  j_i(k)+1} \leq
   \frac{Q_i}{  \Gamma_i(k)^{2(1+\delta)}} \triangleq \alpha_i(k)^2 .$
Thus,  $\left \lfloor \Gamma_i(k)^{2(1+\delta)}\right \rfloor$ steps
   of     \eqref{sa-inner}   {suffice} for obtaining a solution to \eqref{randomized-alg1}  with $ \varepsilon_i(k)$ satisfying
   $ \mathbb{E} \left[    \| \varepsilon_i(k+1)\|^2        \big | \mathcal{F}_k\right]  \leq \alpha_i(k)^2~a.s.$.
   Then  by  the conditional  Jensen's inequality,
    $ \mathbb{E} \left[    \| \varepsilon_i(k+1)\|         \big | 	\mathcal{F}_k\right]  \leq \alpha_i(k)~a.s. $.
        By invoking \eqref{lem-step-size},  we obtain  that for all $i\in \mathcal{N},$
   $\sum_{k=1}^{\infty}\alpha_{i,k }^2<\infty ~a.s.,$ and $\sum_{k=1}^{\infty}\alpha_{i,k }<\infty~a.s.$.
        Then Assumption \ref{assp-noise}(b) holds. \hfill $\Box$
        \end{remark}

 The following result     from  Eqn. (18) in \cite{scutari2010convex}
establishes an    equivalence between  Nash equilibria  of the
stochastic Nash game~\eqref{problem1} and   solutions to the
 variational inequality problem \eqref{optimial-cond}.
  \begin{lemma}\label{lem-equi}~
Let Assumptions  \ref{assump-play-prob}(a), \ref{assump-play-prob}(b),
	and  \ref{assp-potential} hold.    Then    $x^*$ is an NE of the potential game  \eqref{problem1}
if and only if  $x^*$ is a solution to the following  problem:
\begin{equation}\label{optimial-cond}
\nabla_x P(x^*)^T(y-x^*) \geq 0  \quad \forall y \in X.  \end{equation}
Further, the set of Nash equilibria is nonempty and compact.
\end{lemma}

\subsection{Convergence Analysis}\label{Sec:IIC}

We now establish the  almost sure  convergence  and convergence in mean  of the iterates produced by   Alg.  \ref{algo-randomized}.
 Parts of the proof are  inspired by Theorem 4.1 in \cite{davis2016asynchronous} and  Theorem 4.3 in \cite{facchinei2011decomposition}.

 \begin{theorem}[{\bf almost sure  convergence to Nash equilibrium}] \label{thm-rand-0}~
Let  $\{ x(k)\}$ be generated by
 Algorithm \ref{algo-randomized}.  Suppose  Assumptions  \ref{assump-play-prob}, \ref{assp-potential} and  \ref{assp-noise}   hold. We further assume that  {for every $i\in \mathcal{N}$,  the parameter $\mu_i$ utilized  in \eqref{proximal} satisfies $\mu_i>{L_i\over 2}+{ \sqrt{2}   \tau L_i\over 2} ({ L_i \over L_{\rm ave}} +{ L_{\rm ave} \over L_i})$, where $L_{\rm ave}\triangleq \sum_{i\in \mathcal{N}} L_i/N.$} Then the following hold:

\noindent (a) {\bf (square summability):}
For any $i\in \mathcal{N}$,  $  \sum\limits_{k=0}^\infty    \|  T_i(y^i(k) )- x(k) \|^2<\infty \quad a.s.$.

\noindent  (b) {\bf (cluster point is an NE):}
 For almost all $\omega \in \Omega$, every limit point of  $x(k,\omega)$   is a Nash equilibrium.

\noindent  (c)  {\bf (almost sure  convergence to a connected subset of the set $X^*$ of  Nash equilibria):} There exists a connected subset $X_c^* \subset X^*$ such that $d(x(k), X_c^*)  \xlongrightarrow [k \rightarrow \infty]{} 0~a.s.$.
\end{theorem}
{\bf Proof.}  By   Assumption \ref{assump-play-prob}(b),  we have the following bound:
  \begin{equation} \label{fb-term1}
\begin{split}    f_i\left(T_i(y^i(k) ),x_{-i}(k)\right) \leq f_i(x(k))+\nabla_{x_i} f_i( x(k))^T\left(T_i(y^i(k) )-x_i(k)\right)
+ {L_i\over 2} \left\|T_i(y^i(k) )-x_i(k)\right\|^2.
  \end{split}
\end{equation}
 Since  $T_{i}(y^i(k))$ is  a global minimum of \eqref{proximal} and  $x_i(k)\in X_i$,   by the   optimality  condition  we have that
    \begin{align}
0 &\leq     \left(  \nabla_{x_i} f_i( T_{i}(y^i(k)), y_{-i }^i(k))+\mu_i(T_{i}(y^i(k))-x_i(k)) \right)^T\left(x_i(k)-T_{i}(y^i(k))\right)
	\notag \\&= - \left(T_{i}(y^i(k))-x_i(k)\right)^T \nabla_{x_i}   f_i\left(
		T_{i}(y^i(k)), y_{-i}^i(k)\right) -\mu_i\|T_{i}(y^i(k))-x_i(k)\|^2
\notag\\&=- \left(T_{i}(y^i(k))-x_i(k)\right)^T \nabla_{x_i}  f_i( x_i(k),
		y_{-i}^i(k))-\mu_i\|T_{i}(y^i(k))-x_i(k)\|^2  \label{fb-term2}
\\&-\left(\nabla_{x_i}   f_i( T_{i}(y^i(k)), y_{-i}^i(k)) -\nabla_{x_i}
		f_i( x_i(k), y_{-i}^i(k))
		\right)^T\left(T_{i}(y^i(k))-x_i(k)\right) \notag
\\& \leq-\nabla_{x_i}^T  f_i(  y^i(k))\left(T_{i}(y^i(k))-x_i(k)\right)-\mu_i\|T_{i}(y^i(k))-x_i(k)\|^2,\notag
\end{align}
where the last inequality  follows  by  $\big(\nabla_{x_i}   f_i( x_i, x_{-i}) -\nabla_{x_i}   f_i( x_{i}', x_{-i})  \big)^T (x_i-x_i')\geq 0~\forall x_i,x_i'\in X_i,~ \forall x_{-i}\in X_{-i}$ from  Assumption  \ref{assump-play-prob}(b).
    Adding  terms \eqref{fb-term1} and  \eqref{fb-term2},  we have the following inequality:
 \begin{equation} \label{potential-inequ1}
\begin{split}
    f_i\left(T_i(y^i(k) ),x_{-i}(k)\right)  & \leq f_i(x(k))+\left(\nabla_{x_i} f_i( x(k))-\nabla_{x_i}   f_i( y^i(k))\right)^T\left(T_{i}(y^i(k))-x_i(k)\right)\\&-  \left(\mu_i- L_i/2\right) \|T_{i}(y^i(k))-x_i(k)\|^2.
\end{split}
\end{equation}
By Assumption \ref{assump-play-prob}(b) and $ab \leq  {a^2+b^2 \over 2} $, we obtain the following sequences of inequalities for any $C_i>0$:
 \begin{align}
 &  \left(\nabla_{x_i} f_i( x(k))  -\nabla_{x_i}   f_i( y^i(k)) \right)^T\left(T_{i}(y^i(k))-x_i(k)\right) \leq   \left \| \nabla_{x_i} f_i( x(k))  -\nabla_{x_i}   f_i( y^i(k)) \right \| \left\|T_{i}(y^i(k))-x_i(k)\right\| \notag
\\&  \leq   L_i  \|  x(k)-y^i(k)\|\| T_{i}(y^i(k))-x_i(k)\|
 \leq  {L_i^2 \over 2C_i}\|  x(k)-y^i(k)\|^2+  {C_i\over 2}\| T_{i}(y^i(k))-x_i(k)\|^2    . \label{bd-gradient}
\end{align}
 Since   $y^i(k) \triangleq (x_1(k-d_{i1}(k)),\cdots,x_N(k-d_{iN}(k))  )$ and $ d_{ij}(k) \in \{0,1,\cdots,\tau\}  $,  the following holds
\begin{align}
\notag
\|  x(k)-y^i(k)\|^2 &=\sum_{j=1}^N \|  x_j(k)-x_j(k-d_{ij}(k))\|^2= \sum_{j=1}^N \Big \|  \sum_{h=k-d_{ij}(k)+1}^k \left( x_j(h)-x_j(h-1) \right)\Big\|^2
\\& = \sum_{j=1}^N d_{ij}(k)^2  \Big \| {1\over d_{ij}(k)} \sum_{h=k-d_{ij}(k)+1}^k \left( x_j(h)-x_j(h-1) \right)\Big\|^2  \notag
\\& \leq  \sum_{j=1}^N  d_{ij}(k)^2  \Big( {1\over d_{ij}(k)}  \sum_{h=k-d_{ij}(k)+1}^k \|  x_j(h)-x_j(h-1)\|^2 \Big)\quad {\scriptstyle \left(  \mathrm{by ~Jensen's ~inequality}  \right) }\notag
  \\ & \leq  \tau \sum_{j=1}^N    \sum_{h=k-\tau+1}^k \|  x_j(h)-x_j(h-1)\|^2 =   \tau \sum_{h=k-\tau+1}^k \| x(h)-x(h-1)\|^2 \label{seq-inequa}
    \\ \notag &=\tau  \underbrace{    \sum_{h=k-\tau+1}^k (h-k+\tau) \| x(h)-x(h-1)\|^2}_{\triangleq V_k}  +  \tau^2   \| x(k+1)-x(k)\|^2
    \\&\quad -   \tau       \sum_{h=k-\tau+2}^{k+1} \big(h-(k+1)+\tau\big) \| x(h)-x(h-1)\|^2   .\notag
  \end{align}  
 Suppose  $C_i= \sqrt{2} L_i^2\tau /L_{\rm ave}$. Then  by substituting  \eqref{seq-inequa} into \eqref{bd-gradient}, and  by invoking \eqref{potential-inequ1}, we obtain   that
 \begin{equation} \label{potential-inequ2}
\begin{split}   f_i\left(T_i(y^i(k) ),x_{-i}(k)\right)  &\leq f_i(x(k)) +  {L_{\rm ave} \over 2\sqrt{2}} \big(V_k- V_{k+1}+\tau  \| x(k+1)-x(k)\|^2 \big) \\&\quad  - {\left(\mu_i-{L_i+ \sqrt{2} L_i^2\tau /L_{\rm ave} \over 2}\right)}\|T_{i}(y^i(k))-x_i(k)\|^2.
  \end{split}
\end{equation}
  By Algorithm  \ref{algo-randomized},  we have that at time instance $k,$
\begin{equation} \label{diff-triag}
\begin{split}   \| x(k+1)-x(k)\|^2=\|x_{i_k}(k+1)-x_{i_k}(k) \|^2
     \leq 2\|T_{i_k}(y^{i_k}(k))-x_{i_k}(k)\|^2+2\|\varepsilon_{i_k}(k+1)\| ^2.
  \end{split}
\end{equation}
  By     employing  Assumptions  \ref{assump-play-prob}(c), \ref{assump-play-prob}(d),
   and the Jensen's inequality,  the following holds for any $x\in X$.
 \begin{equation}\label{bd-sg}
 \begin{split}
 \|   \nabla_{x_i} f_i(x_i,x_{-i})\| & =\|\mathbb{E}[\nabla_{x_i} \psi_i(x_i, x_{-i};\xi )]\|
   \leq \mathbb{E}[\|\nabla_{x_i} \psi_i(x_i, x_{-i};\xi )\|]
     \\&\leq  \sqrt{\mathbb{E}[\|\nabla_{x_i} \psi_i(x_i, x_{-i};\xi )\|^2]} \leq M.
     \end{split}
       \end{equation}
 Then     by  Algorithm  \ref{algo-randomized} and by invoking   Assumption \ref{assp-potential},  we may obtain
the following bound: \begin{align}
  \ \quad &P(x(k+1) )- P( x(k) )   =P( x_{i_k}(k+1),x_{-i_k}(k)) -  P( x_{i_k}(k),x_{-i_k}(k)) \notag
  \\& =f_{i_k}( x_{i_k}(k+1),x_{-i_k}(k)) -  f_{i_k}( x_{i_k}(k),x_{-i_k}(k)) \notag
  \\&  =   f_{i_k}\left(T_{i_k}(y^{i_k}(k)),x_{-i_k}(k)\right)- f_{i_k}(x(k)) + f_{i_k}( x_{i_k}(k+1),x_{-i_k}(k)) -f_{i_k}\left(T_{i_k}(y^{i_k}(k)),x_{-i_k}(k)\right)
  \label{potential-inequ3}
 \\&    =   f_{i_k}\left(T_{i_k}(y^{i_k}(k)),x_{-i_k}(k)\right)- f_{i_k}(x(k)) +  \varepsilon_{i_k}(k+1)^T \nabla_{x_{i_k}}   f_{i_k} ( z_{i_k}(k+1),x_{-i_k}(k))   \quad {\scriptstyle \left(  \mathrm{by ~mean-value ~theorem}  \right) }  \notag
  \\& \leq   f_{i_k}\left(T_{i_k}(y^{i_k}(k)),x_{-i_k}(k)\right)- f_{i_k}(x(k))+ M \|\varepsilon_{i_k}(k+1)\|,  \quad {\scriptstyle \left(  \mathrm{by ~Cauchy-Schwarz~ inequality~and~} \eqref{bd-sg}  \right) }  , \notag
  \end{align}
       where $z_{i_k}(k+1)=\vartheta_{i_k,k} x_{i_k}(k+1)+(1-\vartheta_{i_k,k}) T_{i_k}(y^{i_k}(k))$ for some $\vartheta_{i_k,k}\in (0,1).$
 Therefore,  by    combining   \eqref{potential-inequ2}, \eqref{potential-inequ3}, and \eqref{diff-triag},  we have the following inequality:
   \begin{equation} \label{potential-inequ4}
\begin{split}  P(x(k+1))  & \leq P(x(k)) + {L_{\rm ave} \over 2\sqrt{2}}(V_k- V_{k+1}) +M \|\varepsilon_{i_k}(k+1)\|  +   {L_{\rm ave}\tau  \over \sqrt{2}}\|\varepsilon_{i_k}(k+1)\| ^2 \\&
- \Big(\mu_{i_k}-{L_{i_k}+\sqrt{2} (L_{i_k}^2/L_{\rm ave}+ L_{\rm ave}) \tau  \over 2} \Big)  \|T_{i_k}(y^{i_k}(k))-x_{i_k}(k)\|^2
  \end{split}
\end{equation}
Therefore, by  rearranging the terms  of \eqref{potential-inequ4} and taking expectations conditioned on  $\mathcal{F}_k $,  we obtain that
\begin{equation}\label{delay-cond-inequ1}
\begin{split}
   & \mathbb{E} \left[P(x(k+1))+ {L_{\rm ave} \over 2\sqrt{2}}V_{k+1}   \big | \mathcal{F}_k \right]
   \leq   \mathbb{E} \left[P(x(k))+ {L_{\rm ave} \over 2\sqrt{2}} V_{k}   \big | \mathcal{F}_k \right]  +   {L_{\rm ave}\tau  \over \sqrt{2}}  \sum_{i=1}^N \mathbb{E} \left[  \|\varepsilon_i(k+1)\|^2     \big | \mathcal{F}_k\right]
\\& + M \sum_{i=1}^N\mathbb{E} \left[   \|\varepsilon_i(k+1)\|      \big | \mathcal{F}_k\right]-     \mathbb{E} \Big[ \Big(\mu_{i_k}-{L_{i_k}+\sqrt{2} (L_{i_k}^2/L_{\rm ave}+ L_{\rm ave}) \tau  \over 2} \Big)
  \|T_{i_k}(y^{i_k}(k))-x_{i_k}(k)\|^2\big | \mathcal{F}_k\Big].
     \end{split}
      \end{equation}
      Since     $T_i (y^i(k))~\forall i\in \mathcal{N}$  is adapted to $\mathcal{F}_k$,
 and  $i_k$ is independent of $\mathcal{F}_k$,
by    Corollary 7.1.2 in \cite{chow2012probability}\footnotemark  \footnotetext{Let the random vectors
$X\in \mathbb{R}^m $ and $Y\in \mathbb{R}^n$ on $(\Omega, \mathcal{F},\mathbb{P})$
be independent of one another and let $f$ be a Borel function on $\mathbb{R}^{m \times n}$ with
$| \mathbb{E}[f(X,Y)] | \leq \infty$.  If for any $x\in  \mathbb{R}^m$, $g(x)=\begin{cases} \mathbb{E}[f(x,Y)]
~&{\rm if~}| \mathbb{E}[f(x,Y)]|\leq \infty \\ 0~&{\rm otherwise} \end{cases}, $
  then $g$ is a Borel function with $g(X)=\mathbb{E}[f(X,Y)| \sigma(X)]$.} and $\mathbb{P}(i_k=i)=p_i$,
 the last term on the right-hand side of \eqref{delay-cond-inequ1} is equivalent to
	    \begin{equation}\label{cond-expect}
\begin{split}
 &  \mathbb{E}_{i_k} \Big[ \Big(\mu_{i_k}-{L_{i_k}+\sqrt{2} (L_{i_k}^2/L_{\rm ave}+ L_{\rm ave}) \tau  \over 2} \Big)   \|T_{i_k}(y^{i_k}(k))-x_{i_k}(k)\|^2  \Big]
\\&=\sum_{i=1}^N p_i  \Big(\mu_{i }-{L_{i}+\sqrt{2} (L_{i}^2/L_{\rm ave}+ L_{\rm ave}) \tau  \over 2} \Big)     \|T_{i}(y^i(k))-x_i(k)\|^2         .
\end{split}
\end{equation}
Since   $ x(k)$ and $V_k$ are adapted to $\mathcal{F}_k$,   by  \eqref{delay-cond-inequ1} and \eqref{cond-expect},  we have the following:
   \begin{equation}     \label{delay-cond-exp}
\begin{split}   & \mathbb{E} \left[P(x(k+1))+ {L_{\rm ave} \over 2\sqrt{2}}V_{k+1}   \big | \mathcal{F}_k \right]
   \leq   \mathbb{E} \left[P(x(k))+ {L_{\rm ave} \over 2\sqrt{2}} V_{k}   \big | \mathcal{F}_k \right] + M \sum_{i=1}^N\mathbb{E} \left[   \|\varepsilon_i(k+1)\|      \big | \mathcal{F}_k\right]
\\& +  {L_{\rm ave}\tau  \over \sqrt{2}} \sum_{i=1}^N  \left[  \|\varepsilon_i(k+1)\|^2     \big | \mathcal{F}_k\right]
-\sum_{i=1}^N p_i  \left(\mu_{i }-{L_{i}+\sqrt{2} (L_{i}^2/L_{\rm ave}+ L_{\rm ave}) \tau  \over 2} \right)     \|T_{i}(y^i(k))-x_i(k)\|^2 .
  \end{split}
\end{equation}

 (a) By using    $\mu_i>{L_i\over 2}+{ \sqrt{2}(L_i^2/L_{\rm ave}+L_{\rm ave} )\tau \over 2}$ and  Assumption \ref{assp-noise}(b),  we may then invoke   Theorem 1  in   \cite{robbins1985convergence},  and conclude that  for every $i\in \mathcal{N}$,  $\sum\limits_{k=0}^\infty   \|T_{i}(y^i(k))-x_i(k)\|^2<\infty,~ a.s . .$

(b) By result (a), we have the following for any $i\in \mathcal{N}$:
 \begin{equation}\label{delay-limit}
 \lim\limits_{k \rightarrow \infty } \|  T_i\left(y^i(k) \right)- x_i(k) \|=0,~ \quad  a.s. \quad .
 \end{equation}
 Let $\bar{x}(\omega)$ be a cluster point of sequence $\{x(k,\omega)\}$.
  Then there   exists a subsequence $\mathcal{K}(\omega)$ such that
 \begin{equation}\label{delay-limit0}
 \lim\limits_{k \rightarrow \infty , k \in \mathcal{K}(\omega)}
x(k,\omega)= \bar{x}(\omega).
\end{equation} Then by \eqref{delay-limit} and \eqref{delay-limit0}, we have that
 \begin{equation} \label{delay-limit1}
 \lim\limits_{k \rightarrow \infty , k \in \mathcal{K}(\omega)} T_i(
		 y^i(k,\omega))= \bar{x}_i(\omega),\quad \forall i\in \mathcal{N}.
\end{equation} We   intend to show that $\bar{x}(\omega)$ is a Nash equilibrium.
  {We proceed by contradication.} Then there exists an $i\in \mathcal{N}$ and a vector $\bar{y}_i \in X_i$ such that
 $ f_i(\bar{y}_i,\bar{x}_{-i}(\omega)) < f_i\left(\bar{x}_i(\omega),\bar{x}_{-i}(\omega)\right) . $
 {By    definition,   the directional derivative of   $f_i $   at point $(\bar{x}_i(\omega),\bar{x}_{-i}(\omega))$ with respect to $x_i$ along the vector $q_i=\bar{y}_i-\bar{x}_i(\omega)$, denoted by $  f_i'(\bar{x}_i(\omega),\bar{x}_{-i}(\omega);q_i)  $, satisfies the following:
 \begin{equation*}
 \begin{split}
  f_i'(\bar{x}_i(\omega),\bar{x}_{-i}(\omega);q_i)   & \triangleq \inf_{\lambda>0} \frac{f_i(\bar{x}_i(\omega)+\lambda q_i,\bar{x}_{-i}(\omega))-f_i(\bar{x}_i(\omega),\bar{x}_{-i}(\omega))}{\lambda}   \\ 	& \leq f_i(\bar{x}_i(\omega)+ q_i,\bar{x}_{-i}(\omega))-f_i(\bar{x}_i(\omega),\bar{x}_{-i}(\omega))
 \quad {\scriptstyle \left(  \mathrm{by ~ setting~}  \lambda=1 \right) }
\\& =f_i(\bar{y}_i,\bar{x}_{-i}(\omega)) - f_i(\bar{x}_i(\omega),\bar{x}_{-i}(\omega))<0 .
\end{split}
\end{equation*}
Since $f_i $ is differentiable at $x$ with respect to $x_i$  by Assumption   \ref{assump-play-prob}(b),  the following holds
 \begin{equation}\label{delay-contradict}
 \begin{split}
 0> f_i'(\bar{x}_i(\omega),\bar{x}_{-i}(\omega);q_i)   &=\left(\bar{y}_i-\bar{x}_i(\omega) \right)^T \nabla_{x_i} f_i(\bar{x}_i(\omega),\bar{x}_{-i}(\omega))     .
\end{split}
\end{equation}}
 Recall that  $T_i(y^i_k)$  is defined as  a global minimum of   a convex optimization problem \eqref{proximal}.
 Since $\bar{y}_i \in X_i$, by  the optimality  condition for a constrained convex programming, we obtain that
   \begin{equation}\label{opt-pbr}
\begin{split}
&\mu_i \left(\bar{y}_i-T_i(y^i(k) ) \right)^T\left(  T_i(y^i(k) ) -x_i(k)\right) +
\left(\bar{y}_i-T_i(y^i(k) ) \right)^T \nabla_{x_i} f_i(T_i(y^i(k) ) ,x_{-i}(k)) \geq 0.
\end{split}
\end{equation}
Since $\nabla_{x_i} f_i(x)$ is continuous in $x$ by Assumption \ref{assump-play-prob}(b),  by taking limits to $k \rightarrow \infty , k \in \mathcal{K}(\omega)$,  and using \eqref{delay-limit}-\eqref{delay-limit1}, we obtain that  $\left(\bar{y}_i- \bar{x}_i (\omega) \right)^T \nabla_{x_i}  f_i(\bar{x}_i(\omega),\bar{x}_{-i}(\omega)) \geq  0,$
   which contradicts \eqref{delay-contradict}.   Thus,   $\bar{x}(\omega)$ is  an NE, proving  (b).

(c)   We  first validate  that $\lim\limits_{k\to \infty} d(x(k), X^*)  = 0~a.s.$.   Assume false; then
$ \limsup  \limits_{k\to \infty} d(x(k,\omega), X^*)  >0 $ with some positive probability,  i.e. for all $\omega \in \widehat{\Omega} \subset \Omega$ where $\mathbb{P}(\widehat{\Omega}) > 0$.  Then   for $\omega \in \widehat{\Omega}$,
by the boundedness of  $\{x(k,\omega)\}_{k\geq 0}$ we can extract a convergent subsequence $\{x(k,\omega) \}_{k\in \mathcal{K}(\omega)}$ such that   $\lim\limits_{k\in \mathcal{K}(\omega), k\to \infty}d(x(k,\omega),  X^*)  >0 $. This contradicts result (b). Hence $\lim\limits_{k\to \infty} d(x(k), X^*)  = 0~a.s.$.

 We now proceed to show a.s. convergence to a  connected  subset of $X^*$, denoted by $X_c^*$, through a contradiction argument.
We assume the contrary   that $X_c^{*}$ is disconnected  with   some positive probability.   Then   there exist at least two  closed connected sets $X_{c_1}^{*}$
 and $X_{c_2}^{*}$ such that $X_{c}^{*}= X_{c_1}^{*} \cup X_{c_2}^{*}$
 with $d(X_{c_1}^{*}, X_{c_2}^{*})>0$.   By hypothesis, the sequence
	 $\{x(k)\}$ cannot converge to either $X_{c_1}^*$ or $X_{c_2}^*$ a.s.,  and $x(k)$ visits $X_{c_1}^{*},X_{c_2}^{*}$ infinitely often.
Define $\rho\triangleq \frac{1}{3}d(X_{c_1}^{*}, X_{c_2}^{*}).$
By $ d( x(k),X_{c}^{*} )  \xlongrightarrow [k \rightarrow \infty]{} 0$,  we know    there exists $k_0$ such that
\begin{equation}
\label{E1}
 x(k)\in B(X_{c_1}^{*}, \rho) \cup   B(X_{c_2}^{*}, \rho) ~~ \forall k\geq k_0,
\end{equation}
where $B(A, \rho)$ denotes the $\rho$-neighborhood   of $A.$
 Define $n_0$ and  $n_p,m_p$ for $p\geq 1$ as follows:
\begin{align*}
& n_0 \triangleq \inf \{ k >k_0, d( x(k),X_{c_1}^{*} )<\rho \},
 m_p \triangleq \inf \{ k >n_{p-1}, d( x(k),X_{c_2}^{*} )<\rho \}, \\
\mbox{ and } & n_p \triangleq \inf \{ k >m_{p}, d( x(k),X_{c_1}^{*} )<\rho \},   ~~ p\geq 1 .
\end{align*}
Then $\{n_p\}$ and $ \{m_p\}$ are infinite sequences by the converse  of result (c).
By \eqref{E1},  we have   $  x(n_p) \in B(X_{c_1}^{*}, \rho)$ and $ x(n_p-1) \in B(X_{c_2}^{*}, \rho)  $
   for any $ p \geq 1$. Then by   $d(X_{c_1}^{*}, X_{c_2}^{*})=3\rho$,   it follows that   $ \|  x(n_p) - x(n_p-1) \| >\rho~\forall p \geq 1 $
   with  some positive probability. Then we have  the following:
   \begin{equation} \label{E2}
 \mathbb{E}[\|  x(n_p) - x(n_p-1) \| ]>0,
\end{equation}
 {where   the  unconditional expectation  is   taken w.r.t.  the   information of $\xi$ and random delays up to time $n_p.$
Hereafter,  the  unconditional expectation   of a variable  is taken w.r.t. to all historical random effects.}

Also, by taking  unconditional  expectations on both sides of \eqref{delay-cond-exp}, we have  the following:
    \begin{equation}     \label{delay-cond-exp1}
\begin{split}    & \mathbb{E} \left[P(x(k+1))+ {L_{\rm ave} \over 2\sqrt{2}}V_{k+1}  \right]
   \leq   \mathbb{E} \left[P(x(k))+ {L_{\rm ave} \over 2\sqrt{2}} V_{k}    \right] + M \sum_{i=1}^N\mathbb{E} \left[   \|\varepsilon_i(k+1)\|     \right]
\\& +    {L_{\rm ave}\tau  \over \sqrt{2}}  \sum_{i=1}^N\mathbb{E} \left[  \|\varepsilon_i(k+1)\|^2     \right]
-\sum_{i=1}^N p_i   \Big(\mu_{i }-{L_{i}+\sqrt{2} (L_{i}^2/L_{\rm ave}+ L_{\rm ave}) \tau  \over 2} \Big)    \mathbb{E} [\|T_{i}(y^i(k))-x_i(k)\|^2] .
  \end{split}
  \end{equation}
   By using Assumption \ref{assp-noise}(b),  we have that
  \begin{equation}   \label{noise-summable1}
\begin{split} &  \sum_{k=1}^{\infty} \mathbb{E}[ \| \varepsilon_i(k+1)\|  ]=\mathbb{E} \Big[\sum_{k=1}^{\infty}  \mathbb{E} \left[    \| \varepsilon_i(k+1)\|        \big | \mathcal{F}_k\right]  \Big]< \infty,   \textrm{~and~} \\&\sum_{k=1}^{\infty} \mathbb{E}[ \| \varepsilon_i(k+1)\|^2  ]=\mathbb{E} \Big[\sum_{k=1}^{\infty}  \mathbb{E} \big[    \| \varepsilon_i(k+1)\|^2        \big | \mathcal{F}_k\Big]  < \infty.
  \end{split}
  \end{equation}
Then from \eqref{delay-cond-exp1} it follows that
 \begin{align*}
& \sum_{k = 0}^\infty\sum_{i=1}^N p_i \Big(\mu_{i }-{L_{i}+\sqrt{2} (L_{i}^2/L_{\rm ave}+ L_{\rm ave}) \tau  \over 2} \Big)   \mathbb{E} [\|T_{i}(y^i(k))-x_i(k)\|^2]    \leq   \mathbb{E}\left[P(x_0) +{L_{\rm ave} \over 2\sqrt{2}}V_0\right]  \\&-
 \liminf_{k \to \infty} \mathbb{E} \left[ P(x(k))+{L_{\rm ave} \over 2\sqrt{2}}V_{k}  \right]   +M \sum_{i=1}^N  \sum_{k = 0}^\infty
    \mathbb{E} \left[    \| \varepsilon_i(k+1)\|  \right]    + {L_{\rm ave}\tau  \over \sqrt{2}} \sum_{i=1}^N  \sum_{k = 0}^\infty
    \mathbb{E} \left[    \| \varepsilon_i(k+1)\|^2       \right]  < \infty,
\end{align*}
where the second inequality holds  by  the   boundedness of $P(x(k))+V_k$ since $P(\cdot)$ is continuous and $X$ is compact.
Hence  by   using $\mu_i>{L_i\over 2}+{ \sqrt{2}(L_i^2/L_{\rm ave}+L_{\rm ave})\tau \over 2}$ and Jensen's inequality,  we have that  \begin{equation}\label{limit-diff0}
 \lim\limits_{k \rightarrow \infty }  \mathbb{E}\left[\|T_{i}(y^i(k))-x_i(k)\|\right]  = 0~~\forall i\in \mathcal{N}.
\end{equation}
 Note that $\| x(k+1)-x(k)\| \leq \sum_{i=1}^N (\| \varepsilon_i(k+1)\|+
		 \|  T_{i}(y^i(k))-x_i(k)\|).$
Then by \eqref{noise-summable1},we have that \begin{equation}\label{limit-diff2}
 \lim\limits_{k \rightarrow \infty } \mathbb{E}[\|x(k+1) -x(k)\|]  = 0.
\end{equation}
This contradicts   \eqref{E2}, and hence the
converse assumption  does not hold. Then result  (c) is proved.
\hfill $\Box$

Theorem \ref{thm-rand-0} shows that  the estimates  generated by Algorithm \ref{algo-randomized}
 converge   almost surely to the set of Nash equilibria.  If the set of Nash equilibria contains isolated points, then for almost all $\omega \in \Omega$,
 $x(k,\omega)$ converges to   an NE. Further, if   the potential game
 \eqref{problem1} admits a unique NE, then  the iterates converge almost surely 
 to the NE.  {It is also worth emphasizing that while we use the term ``best-response'', since $\mu_i$ has to be sufficiently large for all $i$, this can be seen to be more akin to ``better-response.''} In what follows, we  discuss the convergence in mean of the
 iterates.  Since the potential function is employed as a vehicle to
analyze convergence of the iterates, a natural approach would  {have to}
  employ the value of the potential function. However, the
	iterates may converge to stationary points which are not necessarily
	global minimizers and as a consequence, we need to select  an appropriate 	metric to capture stationarity.

Note from Lemma~\ref{lem-equi} that a stationary point of 	$\min_{x \in X} P(x)$ is given by a solution to the variational  inequality problem VI$(X, \nabla_x P)$ that requires an $x \in X$ such that  $ (y-x)^T \nabla_x P(x) \geq 0 ~ \forall y \in X. $
Suppose $X^*$ denotes the set of solutions to VI$(X,F)$. A merit function for ascertaining the departure from solvability of the VI is  a gap function.  It may be recalled from \cite{larsson1994class} that a function $G(x)$ is called a gap
function if it satisfies two properties:
(i) $G(\cdot)$ is sign restricted over the set $X$;   (ii) $G(x) = 0$ if and only if $x$ solves VI$(X,F)$.
We   consider a   primal gap function  \cite[Theorem 3.1]{larsson1994class}  that has found a fair amount of
applicability in the context of variational inequality problems.
\begin{definition} ~Let $X \subseteq \Real^n$ be a nonempty, closed, and convex set.
Let  $F: X \to \mathbb{R}^n $ and let $G:X \to   \mathbb{R}^+$ be defined by $ G(x) \triangleq \sup_{y \in X} F(x)^T(x-y) , \quad \forall x \ \in X.$ 
\end{definition}
The following result shows   the mean convergence in the sense that the limit of  $\mathbb{E}[G(x(k))]$ is zero. This is analogous to showing that expected sub-optimality	  tends to zero in the context of stochastic program.
\begin{theorem}[{\bf Convergence in mean}]\label{thm-rand-1}~
Let  $\{ x(k)\}$ be generated by Algorithm \ref{algo-randomized}.
Suppose    Assumptions  \ref{assump-play-prob}, \ref{assp-potential} and  \ref{assp-noise}   hold, and, in addition, that  {for every $i\in \mathcal{N}$,  the parameter $\mu_i$ utilized  in \eqref{proximal} satisfies $\mu_i>{L_i\over 2}+{ \sqrt{2}(L_i^2/L_{\rm ave}+L_{\rm ave})\tau \over 2}$.}
Then we have that $\lim\limits_{k \to \infty}\mathbb{E} \big[G(x(k))\big]=0.$
\end{theorem}
 {\bf Proof.}  By  Assumption  \ref{assp-potential}, we have the following for any $i\in \mathcal{N}$ and  any $\bar{y}_i\in X_i:$
{\begin{align*}
&\left(x_i(k)-\bar{y}_i  \right)^T  \nabla_{x_i} P(x(k))
=\left(x_i(k)-\bar{y}_i  \right)^T  \nabla_{x_i} f_i(x(k))
   \\&=\left(x_i(k)-T_i(y^i(k) ) \right)^T \nabla_{x_i} f_i(T_i(y^i(k) ) ,x_{-i}(k))
  + \left(T_i(y^i(k) )-\bar{y}_i \right)^T \nabla_{x_i} f_i(T_i(y^i(k) ) ,x_{-i}(k))
  \\& -\left(x_i(k)-\bar{y}_i  \right)^T   \left(\nabla_{x_i} f_i(T_i(y^i(k) ) ,x_{-i}(k))- \nabla_{x_i} f_i(x(k))\right) .
\end{align*}}
Then we have the following sequence of inequalities for any $i\in \mathcal{N}$ and  any $\bar{y}_i\in X_i:$
 \begin{align*}
  &\left(x_i(k)-\bar{y}_i  \right)^T  \nabla_{x_i} P(x(k))  \leq   \mu_i \left(\bar{y}_i-T_i(y^i(k) ) \right)^T\left(  T_i(y^i(k) ) -x_i(k)\right)
\\& +  \left( x_i(k)-T_i(y^i(k) ) \right)^T \nabla_{x_i} f_i(T_i(y^i(k) ) ,x_{-i}(k))
    \\&+\left(\bar{y}_i-x_i(k)  \right)^T   \left(\nabla_{x_i} f_i(T_i(y^i(k) ) ,x_{-i}(k))- \nabla_{x_i} f_i(x(k))\right)     \quad {\scriptstyle (\textrm{by~} \eqref {opt-pbr})}
     \\& \leq    \mu_i \left\| \bar{y}_i-T_i(y^i(k) ) \right\| \left \| T_i(y^i(k) ) -x_i(k)\right\|+  \left \| x_i(k)-T_i(y^i(k) ) \right \| \| \nabla_{x_i} f_i(T_i(y^i(k) ) ,x_{-i}(k))\|   \\&+\left \| \bar{y}_i-x_i(k)  \right \|  \left \| \nabla_{x_i} f_i(T_i(y^i(k) ) ,x_{-i}(k))- \nabla_{x_i} f_i(x(k))\right\|   \quad {\scriptstyle \left(  \mathrm{by ~Cauchy-Schwarz~ inequality} \right) }
      \\& \leq   \left(\mu_i D_{X_i} +M+L_i D_{X_i}\right)  \left \|
		  T_i(y^i(k) ) -x_i(k)\right\|  . \quad {\scriptstyle \left( 	  \mathrm{by  ~Assumptions ~
				  \ref{assump-play-prob}(a),~\ref{assump-play-prob}(b) ~   and~} \eqref{bd-sg}  \right) }.
\end{align*}
 By summing these inequalities over $i$,   it follows that
 {\begin{equation}\label{opt-bound1}
  \begin{split}
 G(x(k))&=\sup_{\bar{y}\in X} \ (x(k)-\bar{y})^T \nabla P(x(k))  =
\sum_{i=1}^N\sup_{\bar{y}_i\in X_i} \ \left(x_i(k)-\bar{y}_i\right)^T  \nabla_{x_i} P(x(k))
 \\& \leq   \sum_{i=1}^N   \left(\mu_i D_{X_i} +M+L_i D_{X_i}\right)  \left \| T_i(y^i(k) ) -x_i(k)\right\|   .
 \end{split}
\end{equation}
  Then, by taking expectations on  both sides of   \eqref{opt-bound1}, we obtain that
 \begin{align*}
      \mathbb{E}\left[  G(x(k)) \right]
     & \leq  \sum_{i=1}^N   \left(\mu_i D_{X_i} +M+L_i D_{X_i}\right) 
	\mathbb{E}\left[\left \| T_i(y^i(k) ) -x_i(k)\right\|  \right]
	\implies
	\lim_{k \to \infty}
      \mathbb{E}\left[  G(x(k)) \right]  \leq 0  \quad {\scriptstyle   \left(
				  \mathrm{by} ~ \eqref{limit-diff0} \right) }.
   \end{align*}
   However, $G(x(k)) \geq 0$ since $x(k) \in X$,  implying  the required result that $ \lim\limits_{k    \to \infty}     \mathbb{E}\left[  G(x(k)) \right]  = 0$.
\hfill $\Box$

We now define an alternative    proximal  gradient-response  map   as follows for $\mu_i>0$:
 \begin{equation}\label{FB}
T_i^{\mu_i}(x)\triangleq \argmin\limits_{y_i  \in X_i}  \left [ (y_i-x_i)^T  \nabla_{x_i} f_i(x) +{\mu_i \over 2} \| y_i-x_i\|^2\right].
\end{equation}
 Since each player's subproblem is convex, definition \eqref{FB} is equivalent to  \begin{equation}\label{FB1}
T_i^{\mu_i}(x)= \Pi_{X_i} \left[   x_i   -{ 1\over \mu_i}   \nabla_{x_i} f_i(x)   \right].
\end{equation}   \vspace{-0.4in}
  { \begin{corollary}[{\bf a.s., mean convergence under proximal gradient-response  map}]
Consider Algorithm \ref{algo-randomized}  with  \eqref{randomized-alg1}  replaced by  the  following variable sample-size projected gradient-response  scheme:
\begin{equation}  \label{PG}
\begin{split}
& x_i(k+1)  = \Pi_{X_i} \Bigg[   x_i(k)  -  {1\over \mu_i}  \Big({1\over {N_i(k)}} \sum_{t=1}^{N_i(k)}\nabla_{x_i}  \psi_i(  x_i(k) ,y^i_{-i}(k);\xi_{i,t}(k))  \Big) \Bigg],
 \end{split}
\end{equation}
where  $N_i(k)\triangleq \left \lfloor  \Gamma_i(k)^{2(1+\delta)}\right \rfloor$   for some positive $\delta>0$ \red{with $\Gamma_i(k)$  defined by \eqref{rem1}}, and  $\xi_{i,1}(k) , \cdots ,\xi_{i,N_i(k)}(k)$  are $N_i(k)$  realizations of the random vector $\xi$. Let Assumptions  \ref{assump-play-prob}, \ref{assp-potential},   \ref{assp-noise}(a),  and  \ref{assp-noise}(c)   hold,
and  $\mu_i>{L_i\over 2}+{ \sqrt{2}(L_i^2/L_{\rm ave}+L_{\rm ave})\tau \over 2}$  for all $i \in \cal N$.
Additionally,   suppose  that for every $i \in \mathcal{N}$, the random variables $\{ \xi_{i,t}(k)\}_{1 \leq t \leq N_i(k)}$ are  independent  of $\mathcal{F}_k$. Then  the results of Theorem  \ref{thm-rand-0}  and Theorem  \ref{thm-rand-1} still hold.
 \end{corollary}}
{\bf Proof.}  {
Define $\varepsilon_i(k+1)\triangleq x_i(k+1)-T_i^{\mu_i}(y^i(k))$. Then by using \eqref{FB1}, \eqref{PG},  and the nonexpansive property of the projection operator  we obtain that $$
\|\varepsilon_i(k+1)\| \leq  {1\over \mu_i {N_i(k)}} \sum_{t=1}^{N_i(k)} \left   \| \nabla_{x_i} \psi_i(  y^i(k);\xi_{i,t}(k))  -\nabla_{x_i} f_i(  y^i(k))  \right \|.$$
Thus, by Assumption   \ref{assump-play-prob}(d), we obtain that
\begin{align*}
 \mathbb{E} \left[    \| \varepsilon_i(k+1)\|^2        \big | \mathcal{F}_k\right]  \leq  {1\over \mu_i^2 N_i(k)}
\mathbb{E}_{\xi}[  \| \nabla_{x_i} \psi_i(  y^i(k);\xi )  -\nabla_{x_i} f_i(  y^i(k))    \|^2]   \leq  {M^2\over \mu_i^2 N_i(k)} ,
\end{align*}
and hence by the conditional Jensen's inequality,  $ \mathbb{E} \left[    \| \varepsilon_i(k+1)\|        \big | \mathcal{F}_k\right]  \leq  {M\over \mu_i \sqrt{N_i(k)}} $.
By noting that $N_i(k)\triangleq \left \lfloor  \Gamma_i(k)^{2(1+\delta)}\right \rfloor$  and  Remark \ref{rem1}(i), Assumption \ref{assp-noise}(b) holds.
 Since $T_i^{\mu_i}(y^i(k))$ is a global minimum of \eqref{FB} and $x_i(k)\in X_i$, by the optimality condition, it follows that

\begin{equation}\label{potential-decrease}
 0\leq-\nabla_{x_i}^T  f_i( y^i(k))\left(T_i^{\mu_i}(y^i(k))-x_i(k)\right)-\mu_i\|T_i^{\mu_i}(y^i(k))-x_i(k)\|^2
\end{equation}
which is indeed the last inequality in Equation \eqref{fb-term2}. Then by inequality \eqref{potential-decrease}, similar to the proof of Theorem  \ref{thm-rand-0}  and Theorem  \ref{thm-rand-1},   we conclude the result.}
\hfill $\Box$

%
%

  \subsection{Generalized Potential   games and Weighted Potential Games}\label{Sec:IID}
We now  consider  the generalized Nash setting  where the strategy sets are coupled in Section \ref{sub-gen},  and   in Section \ref{sec:weight} we consider the weighted potential game, a generalization of standard potential games. 

  \subsubsection{Generalized potential Nash games}\label{sub-gen} We now extend  the separable  constraint
  to the shared constraint regime,  a special case of coupled constraints.  Suppose  there exists  a nonempty closed set $C\in \mathbb{R}^n$ such that  player
  $i$'s feasible set  $X_i(x_{-i})=\{ x_i \in X_i: (x_i,x_{-i})\in C \}$ depends on the rivals' strategies $x_{-i}$, where $X_i \in \mathbb{R}^{n_i}$ are nonempty closed sets   such that $  \prod_{i=1}^N X_i \cap C$ is nonempty.
We say that a point $x\in \mathbb{R}^n$ is feasible if $x_i\in X_i(x_{-i})$ for any $i\in \mathcal{N}.$
The aim of player $i$ is to  choose a strategy $x_i$ that solves the  following parameterized  stochastic program:
 \begin{equation} \label{problem2}
 \min_{x_i  \in X_i(x_{-i})}\quad f_i(x_i,x_{-i})\triangleq  { \mathbb{E}_{\xi}}\left[
\psi_i(x_i,x_{-i};\xi(\omega)) \right].
\end{equation}
Assume that there exists a continuous  potential function  $P: C\bigcap   \prod_{i=1}^N   X_i \rightarrow \mathbb{R}$
such that for  any $i\in \mathcal{N}$ and any $x_{-i} $ with  $X_i(x_{-i})$ being nonempty, we have the following equality:
\begin{align}\label{genr-potential}  P(x_i,x_{-i}) - P(x_i',x_{-i}) &= f_i(x_i,x_{-i})
-f_i(x_i',x_{-i}), ~~\forall x_i, x_i' \in X_i(x_{-i}).
\end{align}
 Then the problem \eqref{problem2}  is called a generalized    potential  stochastic game where the (generalized) NE $x^* $   is a feasible point  such that  the following holds  for any $i$:
\begin{align*}  	f_i(x_i^*,x_{-i}^*) \leq f_i(x_i ,x_{-i}^*)\quad \forall x_i  \in X_i(x_{-i}^*).
\end{align*}

Suppose  that for any $i\in \mathcal{N}$,  {(i) $f_i(x)$ is continuously  differentiable   on $C$,
(ii)   the feasible set $X_i(x_{-i})$ is inner-semicontinous relative to the set of points $x_{-i}$ for which  $X_i(x_{-i})$ is nonempty
\cite[A2]{facchinei2011decomposition},  (iii)  $X_i(x_{-i})$  is convex for all  $x_{-i}$ for which  $X_i(x_{-i})$ is nonempty, and $f_i(\cdot,x_{-i})$  is  convex in  $x_i \in X_i(x_{-i})$ for  all  $x_{-i}$ for which  $X_i(x_{-i})$ is nonempty}.
 Then computing  the proximal  best-response  solution $T_i(x)$ of   \eqref{problem2}  involves
 solving a strongly convex stochastic  program.
 Let Algorithm  \ref{algo-randomized} with  $d_{ij}(k)=0~\forall i,j\in \mathcal{N}$   be applied to the generalized stochastic potential game with   feasible initial point $x_0$, where it is required that   $x_i(k+1)\in X_i(x_{-i}(k))$,  which can be guaranteed by using the   projected stochastic gradient scheme to  obtain the  approximate proximal BR solutions.
Similar to Lemma 4.1 in \cite{facchinei2011decomposition}, we may  show  that
 $x(k+1)$   is feasible. We then conclude the following result, for which the proof is similar  to that of Theorems  \ref{thm-rand-0}  and   \ref{thm-rand-1}.

 \begin{corollary}[{\bf Generalized potential  stochastic Nash games}]
Let Algorithm  \ref{algo-randomized}   be applied to   the stochastic generalized  Nash game \eqref{problem2}  satisfying \eqref{genr-potential},  where  $x_0$ is feasible,  $y^i(k)=x(k)  $, and  it is  required that   $x_i(k+1)\in X_i(x_{-i}(k))$.
Then the results of Theorems  \ref{thm-rand-0}  and   \ref{thm-rand-1} still hold  under suitable conditions.
 \end{corollary}

   {Nevertheless,  we may be unable to manage the delayed regime  
since we   cannot guarantee that  the delay-afflicted
rival strategies,  denoted by $y^i_{-i}(k)$, allow for retaining feasibility; 
namely,   for some $i, k$,  the set $X_i(y_{-i}^i(k) )$  may be empty.  Consequently,  player $i$  may be  unable to find a feasible strategy, given rival
strategies, and hence  Algorithm  \ref{algo-randomized}    is
 not well-defined in that  step (S.2) cannot be implemented. }

  \subsubsection{Weighted potential games} \label{sec:weight} We now consider the weighted  potential game,
in which there exist positive numbers $w_1,\cdots, w_N$  such that, for  any $i\in \mathcal{N}$ and any $x_{-i} $ the following equality holds:
\begin{align} \label{weighted-potential}  P(x_i,x_{-i}) - P(x_i',x_{-i}) &=w_i\left( f_i(x_i,x_{-i})-f_i(x_i',x_{-i}) \right) ~~\forall x_i, x_i' \in X_i(x_{-i}).
\end{align}
Then by applying Algorithm  \ref{algo-randomized},  we obtain  the following  results.

 \begin{corollary}[{\bf Weighted potential stochastic Nash games}]~
Let Algorithm  \ref{algo-randomized}  be applied to the stochastic  Nash game \eqref{problem1}, in which   the objective function
satisfies   \eqref{weighted-potential}.  Suppose    Assumptions  \ref{assump-play-prob}, \ref{assp-potential} and  \ref{assp-noise}   hold.
  {Let the parameter $\mu_i$ used  in \eqref{proximal} satisfy  $\mu_i>{L_i\over 2}+{\sqrt{2}\tau L_i\over 2} (
{L_{\rm ave}  \over  L_i}{  w_{\rm ave}\over  w_i }+ {  L_i    \over L_{\rm ave}}{w_i \over w_{\rm ave}})$, where $w_{ \rm ave}=\sum_{i\in \mathcal{N}} w_i/N$.}
Then the   results of  Theorem \ref{thm-rand-0}  and Theorem  \ref{thm-rand-1}  hold.
 \end{corollary}
 {\bf Proof.}
 By multiplying both sides of \eqref{potential-inequ1} with $w_i$ and rearranging the terms  we obtain that
 \begin{align}&  w_i\Big(  f_i\left(T_i(y^i(k) ),x_{-i}(k)\right)  - f_i(x(k)) \Big) \notag \\
 &\leq w_i\left(\nabla_{x_i} f_i( x(k))-\nabla_{x_i}   f_i( y^i(k))\right)^T\left(T_{i}(y^i(k))-x_i(k)\right)-w_i
 (\mu_i-L_i/2) \|T_{i}(y^i(k))-x_i(k)\|^2 \label{potential-inequ-w}.
 \end{align}
%
By substituting \eqref{bd-gradient} and \eqref{seq-inequa} into \eqref{potential-inequ-w},
 the following  holds for  $C_i={\sqrt{2} L_i^2 w_i \tau \over L_{\rm ave} w_{\rm ave}}$: 
\begin{equation} \label{potential-inequ-w1}
\begin{split}
  &   w_i\Big(  f_i\left(T_i(y^i(k) ),x_{-i}(k)\right)  - f_i(x(k)) \Big)\leq
  {w_iL_i^2 \tau \over 2C_i}  (V_k-V_{k+1} +  \tau    \| x(k+1)-x(k)\|^2)  \\& -w_i
 {(\mu_i-L_i/2)}\|T_{i}(y^i(k))-x_i(k)\|^2+   {w_iC_i\over 2}\| T_{i}(y^i(k))-x_i(k)\|^2.
\end{split}
\end{equation}
 Then by \eqref{weighted-potential}, similar to  \eqref{potential-inequ3} we have that
 \begin{equation*}
\begin{split}
  &P(x(k+1) )- P( x(k) )    \leq  w_{i_k}  \big( f_{i_k}\left(T_{i_k}(y^{i_k}(k)),x_{-i_k}(k)\right)- f_{i_k}(x(k))+ M \|\varepsilon_{i_k}(k+1)\| \big),
  \end{split}\end{equation*}
which incorporating  with  \eqref{diff-triag} and   \eqref{potential-inequ-w1} yields  the following inequality:
   \begin{equation} \label{potential-inequ-w2}
\begin{split}
 P(x(k+1))-P(x(k))&\leq {L_{\rm ave} w_{\rm ave} \over 2\sqrt{2}} (V_k-V_{k+1})+{L_{\rm ave} w_{\rm ave} \tau \over \sqrt{2}} \|\varepsilon_{i_k}(k+1)\|^2 + M w_{i_k}\|\varepsilon_{i_k}(k+1)\|
\\&-  w_i \left((\mu_i- L_i / 2 )-{L_{\rm ave} w_{\rm ave} \tau\over \sqrt{2}w_i }- { L_i^2 w_i  \tau \over \sqrt{2} L_{\rm ave} w_{\rm ave}}\right)  \|T_{i_k}(y^{i_k}(k))-x_{i_k}(k)\|^2.
\end{split}
\end{equation}   Since $\mu_i>{L_i\over 2}+{\tau L_i\over 2} (
{L_{\rm ave}  w_{\rm ave} \over  L_i   w_i }+ {  L_i     w_i \over  L_{\rm ave} w_{\rm ave}})$,  by using   \eqref{potential-inequ-w2},  similar to the proof of Theorems \ref{thm-rand-0} and  \ref{thm-rand-1},  we  obtain  the results.
\hfill $\Box$

 {
\subsection{ Asynchronous inexact  best-response  scheme  without delays}\label{sec:Br}
We now  show that the asynchronous inexact  pure  best-response  scheme  (i.e.  without a proximal term)  is also applicable when
 the player-specific problem  is strongly convex and  each player may obtain its rivals' latest strategies without  delay.
We define   best-response map  of player $i$ as follows:
\begin{align}\label{def-br}
 \widehat{T}_i(x) \triangleq \argmin_{y_i \in X_i}  \left[ f_i(y_i, x_{-i})  \right],
\end{align}
 based  on which,  analogous  to Algorithm \ref{algo-randomized},  we  design an asynchronous inexact best-response  scheme. }
\vspace{-0.2in}

 \begin{algorithm} [H]
\caption{ Asynchronous inexact best-response scheme }  \label{alg2}
 Let $k:=0$, $ x_{i,0} \in X_i$  for $ i\in \mathcal{N}$.
		 Additionally $0 < p_i < 1$ for $ i\in \mathcal{N}$ such that
$\sum_{i=1}^N  p_i = 1.$
\begin{enumerate}
\item[(S.1)] Pick $i_k=i \in \mathcal{N} $ with probability $p_i$.
\item[(S.2)] If $i_k=i$,  then    player $i$  updates     $ x_i(k+1) \in X_i$  as follows:
\begin{equation*}
x_i(k+1) :=  \widehat{T}_i(x(k))+ \varepsilon_i(k+1),
\end{equation*} where  $  \varepsilon_i(k+1) $ denotes the inexactness.
 Otherwise,     $ x_j(k+1):= x_j(k)$ if $j  \notin  i_k$.
 \item[(S.3)]  {If $k > K$, stop; Else, $k:=k+1$} and return to (S.1).
\end{enumerate}
\end{algorithm}
\vspace{-0.2in}
  {
 \begin{theorem}\label{thm-br}
Let  $\{ x(k)\}$ be generated by  Algorithm \ref{alg2}.  Suppose  Assumptions  \ref{assump-play-prob}, \ref{assp-potential}, and \ref{assp-noise}(b)   hold.
Assume that  for every $i \in \mathcal{N},$ $\sum_{k=1}^{\infty}  \mathbb{E} \left[    \| \varepsilon_i(k+1)\|        \big | \mathcal{F}_k\right]  < \infty~ a.s., $  and    $f_i(x_i,x_{-i})$ is  $\mu_i$-strongly  convex    in $x_i \in X_i$ for  all $x_{-i} \in X_{-i}$. Then for almost all $\omega \in \Omega$, every limit point of  $\{x(k,\omega)\}$   is an NE.
\end{theorem}
}
{\bf Proof.} 
 {
Since  $\widehat{T}_{i}(x(k))$ is  a global minimum  of the problem \eqref{def-br},   by  $x_i(k)\in X_i$ and  the   optimality  condition,  we obtain that $     \big(x_i(k)-\widehat{T}_{i}(x(k))\big)^T \nabla_{x_i} f_i(\widehat{T}_{i}(x(k)),x_{-i }(k)) \geq 0. $  Then by the $\mu_i$-strong   convexity of $f_i(x_i,x_{-i})$ in $x_i \in X_i$, we have that
 \begin{align*} f_i(x(k)) &\geq f_i\left(\widehat{T}_{i}(x(k)),x_{-i}(k)\right) +\left(x_i(k)-\widehat{T}_{i}(x(k))\right)^T\nabla_{x_i}  f_i\left(\widehat{T}_{i}(x(k)),x_{-i}(k)\right)
+ {\mu_i\over 2} \left\|\widehat{T}_{i}(x(k))-x_i(k)\right\|^2
\\&  \geq f_i\left(\widehat{T}_{i}(x(k)),x_{-i}(k)\right)
+ {\mu_i\over 2} \left\|\widehat{T}_{i}(x(k))-x_i(k)\right\|^2.
\end{align*}
By rearranging the terms and using  Assumption \ref{assump-play-prob}(b), we have the following:
\begin{align*}
 &  f_i\left(\widehat{T}_{i}(x(k)),x_{-i}(k)\right)  \leq   f_i(x(k))  -   {\mu_i\over 2} \left\|\widehat{T}_{i}(x(k))-x_i(k)\right\|^2.
\end{align*}
Similar to \eqref{potential-inequ3}, we can also show that
\begin{align*}
  \ \quad &P(x(k+1) )- P( x(k) )     =f_{i_k}( x_{i_k}(k+1),x_{-i_k}(k)) -  f_{i_k}( x_{i_k}(k),x_{-i_k}(k))
  \\& \leq   f_{i_k}\left(\widehat{T}_{i_k}(x(k)),x_{-i_k}(k)\right)- f_{i_k}(x(k))+ M \|\varepsilon_{i_k}(k+1)\|
\\&  \leq  -   {\mu_{i_k}\over 2} \left\|\widehat{T}_{i_k}(x(k))-x_{i_k}(k)\right\|^2+ M \|\varepsilon_{i_k}(k+1)\|.
  \end{align*}
By taking expectations conditioned on  $\mathcal{F}_k $, similarly to \eqref{delay-cond-exp}, we conclude that
\begin{align*}
  & \mathbb{E} \left[P(x(k+1))  \big | \mathcal{F}_k \right]
   \leq  P(x(k))  + M \sum_{i=1}^N\mathbb{E} \left[   \|\varepsilon_i(k+1)\|      \big | \mathcal{F}_k\right]
-\sum_{i=1}^N {p_i \mu_i \over 2} \left\|\widehat{T}_{i}(x(k))-x_i(k)\right\|^2  .
\end{align*}
Since $\sum_{k=1}^{\infty}  \mathbb{E} \left[    \| \varepsilon_i(k+1)\|        \big | \mathcal{F}_k\right]  < \infty ~a.s.  $ for every $i\in \mathcal{N}$,  we   then use  Theorem 1  in   \cite{robbins1985convergence},  and conclude that  $
\sum\limits_{k=0}^\infty   \big\|\widehat{T}_{i}(x(k))-x_i(k)\big\|^2<\infty ~ a.s .$ for every $i\in \mathcal{N}$.
Thus,  $\lim\limits_{k\to \infty}\|\widehat{T}_{i}(x(k))-x_i(k) \|=0 ~ a.s . $ for every $i\in \mathcal{N}$.
The rest of the proof is the same as  that of Theorem \ref{thm-rand-0}(b).
}
\hfill $\Box$

 {
\begin{remark} Theorem \ref{thm-br} shows that in delay-free regimes,   asynchronous  inexact  pure  BR  schemes    retain almost sure  convergence   when player-specific problems are strongly convex.  Consequently, in deterministic regimes,  the exact BR scheme is  convergent when each player's subproblem is strongly convex, complementing the findings from \cite{facchinei2011decomposition} in that the best-response schemes can  lead to  convergence to Nash equilibria when the player-specific problem is   strongly convex rather than merely convex.
\end{remark}}

	\section{Misspecified  Potential  Stochastic Nash Games}\label{sec:misp}	
In this section,  we consider  the misspecified stochastic Nash
game \textbf{(P2)}.  {A sequential approach for resolving such a problem relies on first estimating $\theta^*$ and subsequently estimating $x^*$ based on the belief regarding $\theta^*$.  As pointed by
\cite{ahmadi2014resolution,jiang2017distributed},   this sequential
approach   is characterized by several shortcomings: (i)  In any sequential
approach, the computation of  $\theta^*$ has to be completed in  finite time;
this is generally  impossible  since  $\theta^*$  is defined as a solution to
the stochastic program.  (ii) If the learning of $\theta^*$ is terminated
prematurely, this leads to  an erroneous estimate $\hat{\theta}$.  One then
proceeds to compute a Nash equilibrium  given $\hat{\theta}$, which  results
in an incorrect Nash equilibrium.  As a result, the two-stage sequential
method,  in stochastic regimes,   cannot provide asymptotically accurate
solutions and at best provided approximate solutions.  Motivated by these
shortcomings,  we propose a  framework that combines the  asynchronous
inexact proximal best-response  scheme with joint  learning for the
misspecified  parameter $\theta^*$.  Under suitable conditions, we prove the almost sure convergence and  the  convergence in mean of the generated  strategy
vector to the set of Nash equilibria.  Additionally, we  show that for every
player,  its belief   regarding the misspecified parameter  converges almost surely   to the true counterpart.  }

\subsection{ Algorithm Design  and Assumptions}    We impose the following conditions on the misspecified  problem.
	\begin{assumption} \label{assp-miss}~
	(a) For every $i\in \mathcal{N},$  $X_i $ is a  closed,  compact, and    convex set;
	 $f_i(x_i,x_{-i};\theta)$ is convex   and continuously  differentiable
	  in $x_i$   over an open set containing $X_i$ for every $x_{-i} \in X_{-i}$ and every $\theta\in \Theta.$

\noindent  (b)   For    every $i\in \mathcal{N},$ $\nabla_{x_i} f_i(x;\theta^*)$ is Lipschitz continuous
in $x $ with Lipschitz constant $L_x$, i.e.,
$$\| \nabla_{x_i} f_i(x;\theta^*)-\nabla_{x_i} f_i(x';\theta^*)\| \leq L_x\| x-x'\|\quad \forall x,x'\in X .$$
 Further,  there exists a  constant  $L_{\theta^*}$ such that for any $x\in X$ and every $i\in \mathcal{N}:$
$$\| \nabla_{x_i} f_i(x_i,x_{-i};\theta)-\nabla_{x_i} f_i(x_i,x_{-i};\theta^*)\| \leq  L_{\theta^*} \| \theta-\theta^*\|\quad \forall \theta \in \Theta.$$

\noindent  (c)   $g(\theta)$ is strongly  convex  with convexity constant $\mu_g$   and
is  continuously  differentiable in    $\theta$ \us{on an open set containing} $\Theta $ with the gradient function being   $L_g$-Lipschitz continuous,  {where $g(\theta)$ is defined in \eqref{g_theta}.}

\noindent (d) There exists  a  function
 $P(\cdot;\cdot ):X \times \Theta \to \mathbb{R}$  	 such that for every $ i\in \mathcal{N}$ and any  $x_{-i} \in X_{-i}$,
\eqref{def-potential} holds.
\end{assumption}
 We  define $ P(x)\triangleq  P(x;\theta^*)$ as the potential function of the problem \eqref{problem1-mis}.

\begin{assumption} \label{assp-diff}~
	(a)  For any $ i\in \mathcal{N}$,  all $x_{-i} \in X_{-i}  $, any $\theta\in \Theta$ and  any $\xi \in  \mathbb{R}^d$, $  \psi_i(x_i, x_{-i};\theta;\xi)$ is differentiable  in $x_i$ over {an open set containing} $X_i$ such that
$ \nabla_{x_i} f_i(x_i,x_{-i};\theta)=  \mathbb{E}[\nabla_{x_i} \psi_i(x_i, x_{-i};\theta;\xi  )].$
(b) For   any $ i\in \mathcal{N}$ and  any  $x \in X $,  there exists a constant $M_1>0$ such that
 $\mathbb{E}[\|\nabla_{x_i} \psi_i(x_i,x_{-i};\theta;\xi)\|^2] \leq  M_1^2.$
(c)  For    any $\eta \in \mathbb{R}^p$, $ g(\theta,\eta )$ is differentiable  in $\theta$ over
an open set containing $\Theta$ such that
$ \nabla g(\theta)=  { \mathbb{E}_{\eta}}[\nabla g(\theta;\eta )].$
\end{assumption}

	 If $T_i(x,\theta)$ is defined as follows:
	\begin{equation}\label{def-br2}T_i(x,\theta) \triangleq \argmin_{y_i \in X_i}  \left[ f_i(y_i, x_{-i};\theta ) +\frac{\mu}{2} \|y_i-x_i\|^2\right],\quad \mu>0,
	\end{equation}
then $T_i(x,\theta)$ is uniquely defined  by  invoking Assumption \ref{assp-miss}(a).
 Additionally, we may claim} the Lipschitz continuity   of  $T_i(x,\cdot)$ \us{based}  on   the following Lemma, akin to the result proved by ~\cite{dafermos1988sensitivity}. 

 {\begin{lemma}\label{Lip-Tmap}~
Define   $L_t \triangleq {\mu  L_{\theta^*} \over \mu^2+L_x^2}\big(1-  {L_x  / \sqrt{ \mu^2+ L_x^2} } \big)^{-1}  $. Then  for any $i\in \mathcal{N}$ and any $x\in X$:
	\begin{align}\label{T-Lipschitz}
\|T_i(x,\theta)-T_i(x,\theta^*)\|\leq L_t\|  \theta-\theta^* \|\quad \forall \theta\in \Theta.
\end{align}
\end{lemma} } 
 {\bf Proof.}
By the first-order optimality condition of \eqref{def-br2},  $T_i(x,\theta)$  is  a fixed point of  the map $ \Pi_{X_i} \big[y_i-\alpha \left( \nabla_{x_i}    f_i(y_i, x_{-i};\theta ) + \mu  (y_i-x_i) \right ) \big]$. Then  by  using  the nonexpansivity property of the projection operator, the triangle  inequality,  and  Assumption \ref{assp-miss}(b), we have that
\begin{equation}\label{dis-t}
\begin{split}
&\|T_i(x,\theta)-T_i(x,\theta^*)\|   =\Big\| \Pi_{X_i} \big[T_i(x,\theta)-\alpha \left( \nabla_{x_i}    f_i(T_i(x,\theta), x_{-i};\theta ) + \mu  (T_i(x,\theta)-x_i) \right ) \big]
\\&\qquad- \Pi_{X_i} \big[T_i(x,\theta^*)-\alpha  \left( \nabla_{x_i}   f_i(T_i(x,\theta^*), x_{-i};\theta^* ) + \mu (T_i(x,\theta^*)-x_i) \right) \big] \Big\|
\\&\leq \Big\|(1-\alpha \mu) \left(T_i(x,\theta)-T_i(x,\theta^*) \right)-\alpha\left(  \nabla_{x_i}   f_i(T_i(x,\theta), x_{-i};\theta^* ) - \nabla_{x_i}   f_i(T_i(x,\theta^*), x_{-i};\theta^*  )  \right)
\\& \qquad-\alpha\left( \nabla_{x_i}   f_i(T_i(x,\theta), x_{-i};\theta)- \nabla_{x_i}   f_i(T_i(x,\theta ), x_{-i};\theta^*  ) \right)  \Big\|
 \leq \alpha  L_{\theta^*} \| \theta-\theta^*\|  \\
& + \Big\|(1-\alpha \mu) \left(T_i(x,\theta)-T_i(x,\theta^*) \right)-\alpha\left(  \nabla_{x_i}   f_i(T_i(x,\theta), x_{-i};\theta^* ) - \nabla_{x_i}   f_i(T_i(x,\theta^*), x_{-i};\theta^*  )  \right) \Big\| .
\end{split}
\end{equation}
By recalling that   $f_i(x_i,x_{-i};\theta)$ is convex     $x_i \in X_i$ for all  $x_{-i} \in X_{-i}$, $\theta\in \Theta,$  and  $\nabla_{x_i} f_i(x;\theta^*)$ is Lipschitz continuous in $x $ with Lipschitz constant $L_x$,
  we conclude that for    $\alpha={\mu \over \mu^2+L_x^2}$,
\begin{align*}
&  \Big\|(1-\alpha \mu) \left(T_i(x,\theta)-T_i(x,\theta^*) \right)-\alpha\left(  \nabla_{x_i}   f_i(T_i(x,\theta), x_{-i};\theta^* ) - \nabla_{x_i}   f_i(T_i(x,\theta^*), x_{-i};\theta^*  )  \right) \Big\|^2
\\& \leq \left \| (1-\alpha \mu) (T_i(x,\theta)-T_i(x,\theta^*) \right\|^2+ \alpha^2 \left \|  \nabla_{x_i}   f_i(T_i(x,\theta), x_{-i};\theta^* ) - \nabla_{x_i}   f_i(T_i(x,\theta^*), x_{-i};\theta^*  ) \right\|^2
\\& -2\alpha (1-\alpha \mu) \left(T_i(x,\theta)-T_i(x,\theta^*) \right)^T\left(  \nabla_{x_i}   f_i(T_i(x,\theta), x_{-i};\theta ^*) - \nabla_{x_i}   f_i(T_i(x,\theta^*), x_{-i};\theta^*  )  \right)
\\& \leq (  (1-\alpha \mu)^2+\alpha^2L_x^2) \left \| T_i(x,\theta)-T_i(x,\theta^*)\right\|^2=
 {L_x^2 \over \mu^2+ L_x^2}  \left \|  T_i(x,\theta)-T_i(x,\theta^*) \right\|^2.
\end{align*}
This   incorporated  with  \eqref{dis-t}   implies that  the following holds for $\alpha={\mu \over \mu^2+L_x^2}$:
$$\Big(1-  {L_x  \over \sqrt{ \mu^2+ L_x^2} } \Big)\|T_i(x,\theta)-T_i(x,\theta^*)\| \leq {\mu  L_{\theta^*} \over \mu^2+L_x^2}  \| \theta-\theta^*\|.$$
The  result follows   by the definition of $L_t.$
\hfill $\Box$

	We  propose an   asynchronous   inexact    proximal BR    scheme that is coupled  with
learning  (Algrithm \ref{algo-randomized2})  to   compute  an NE of the misspecified potential stochastic  game.
 {Player $i$  at time $k$  utilizes  an estimate $x_i(k) $ of  its equilibrium strategy $x_i^* $, an estimate $\theta_i(k)$ of  the unknown parameter $\theta^*,$  and has access to   possibly delay-afflicted   rival strategies  $y^i(k)\triangleq (x_1(k-d_{i1}(k)),\cdots,x_N(k-d_{iN}(k)) )$ with delays $d_{ij}(k),j\in\mathcal{N}.$     The scheme  is
defined as follows:   At major iteration $k\geq 0,$   he index $i$ is selected randomly from  $
\mathcal{N} $ with  probability   $\mathbb{P}(i_k=i)=p_i>0$.   If $i_k=i,$  then player $i$ is
chosen to  initiate  an update  by   computing an inexact
proximal  BR  solution to  the problem \eqref{def-br2}  characterized by
\eqref{randomized2}, and   updating   $\theta_i(k+1)$   via the variable
sample-size SA scheme  \eqref{VSSG} with $N_i(k)$ sampled  gradients.
 We    impose conditions on the inexactness sequence    $ \{
\varepsilon_i(k)\}_{k\geq 1} $ and  further  specify the
selection of     $N_i(k) $ in the convergence analysis.}

 \begin{algorithm}  [H]
\caption{{Asynchronous   inexact proximal best-response scheme with stochastic learning }}\label{algo-randomized2}
 Let $k:=0$, \red{$ x_i(0) \in X_i$ and $ \theta_i(0) \in \Theta$}  for $ i\in \mathcal{N}$.   Additionally $0 < p_i < 1$ for $ i\in \mathcal{N}$ such that $\sum_{i=1}^N  p_i = 1.$
		  \begin{enumerate}
\item[(S.1)] Pick $i_k=i \in \mathcal{N} $ with probability $p_i$.
\item[(S.2)] If $i_k=i$,  then    player $i$  updates     $ x_i(k+1) \in X_i $ and $\theta_i(k+1) \in \Theta $  as follows:
\begin{align}
x_i(k+1) &:=  T_i(y^i(k),\theta_i(k))+ \varepsilon_i(k+1),  \label{randomized2}
\\ \theta_i(k+1) &:=\Pi_{\Theta}\Big[\theta_i(k)-\frac{\beta_i}{N_i(k)} \sum_{p=1}^{N_i(k)} \nabla  g\left(\theta_i(k),\eta_{i,p}(k) \right)\Big], \label{VSSG}
\end{align}
where  $  \varepsilon_i(k+1) $ denotes the inexactness,   $  \nabla g\left(\theta_i(k),\eta_{i,p}(k) \right), p=1,\cdots, N_i(k)$ denotes the sampled gradient, \red{and $\beta_i={1\over L_g}$};
 Otherwise,     $ x_j(k+1):= x_j(k),  \theta_{j,k+1}= \theta_j(k)$ if $j  \neq  i_k$.
 \item[(S.3)]  If $k > K$, stop; Else, $k:=k+1$ and return to (S.1).
\end{enumerate}
\end{algorithm}

We  then  list the following conditions concerning  the  delays,  observation noise of the
gradient function $\nabla g(\theta) $  as well as
the inexactness sequence  $\{ \varepsilon_i(k)\}$ utilized  in Algorithm   \ref{algo-randomized2}. {We denote   the $\sigma$-field of the entire information  used by Algorithm \ref{algo-randomized2} up to  (and including)  the updates of  $x_i(k),\theta_i(k) $ for all $ i\in \mathcal{N}$ by $\mathcal{F}'_k$, and  the   $\sigma$-field  generated from $\mathcal{F}'_k$  and the     delays  at step $k$ by   $\mathcal{F}_k\triangleq \sigma \big\{ \mathcal{F}_k',    d_{ij}(k), i,j\in \mathcal{ N}  \big\}$. We will  further  define  $\mathcal{F}'_k$ after introducing the SA scheme \eqref{SA2}.}
 \begin{assumption} \label{assp-noise2}

(a) $\{i_k\}$ is an i.i.d. sequence,   where  $i_k$ is independent of $\mathcal{F}_k$ for all $k\geq 1.$

\noindent (b) For any $i \in \mathcal{N},$ the noise term $\{ \varepsilon_i(k)\}$
 satisfies the following condition:
 $$ \sum_{k=1}^{\infty}  \mathbb{E} \left[    \| \varepsilon_i(k+1)\|^2        \big | \mathcal{F}_k\right]  < \infty,
~ \textrm{and}~\sum_{k=1}^{\infty}  \mathbb{E} \left[    \| \varepsilon_i(k+1)\|        \big | \mathcal{F}_k\right]  < \infty~~a.s. .  $$

\noindent    (c)   There exists a positive integer $\tau$
such that for any $i,j\in \mathcal{N}$ and any $k\geq0,$ $d_{ij}(k) \in \{0,\cdots, \tau\}$.

\noindent  (d) Define $e_{i,p}(k) \triangleq  \nabla g\left(\theta_i(k),\eta_{i,p}(k) \right)- \nabla g\left(\theta_i(k)\right)$.
 There exists a constant $M_2>0$ such that \red{for any $k\geq 0$ and $p=1,\cdots,N_i(k),$}
 $ \mathbb{E} \left[\|e_{i,p}(k)  \|^2| \mathcal{F}_k\right] \leq  M_2^2  .$ 
\end{assumption}

Analogous  to  the computation of   the  inexact best-response   \eqref{randomized-alg1}  in Algorithm \ref{algo-randomized}, we still utilize SA
to   compute \eqref{randomized2}.   By \eqref{def-br2}  it is seen   that  the computation of  $T_i(x,\theta) $  requires solving a strongly convex stochastic program.
Thus, an  inexact solution  to  the problem  \eqref{def-br2}, characterized by
\eqref{randomized2}, can also be computed via the  SA algorithm defined  as follows    for   $t = 1, \hdots, j_i(k)$: 
\begin{equation} \label{SA2}
\begin{split}
& x_{i,t+1}(k)  := \Pi_{X_i} \Big[  z_{i,t}(k) -   \gamma_{i,t}  \left[\nabla_{x_i}
 \psi_i(  z_{i,t}(k) ,y^i_{-i,k};\theta_i(k);\xi_{i,t}(k) ) + \mu (  z_{i,t}(k)-x_i(k))\right] \Big],
 \end{split}
\end{equation}
 where  $ \gamma_{i,t}=\frac{1}{\mu (t+1)}$, and  { $ z_{i ,t}(k) $   denotes  the  estimate  of the proximal BR  solution  $ T_i(y^i(k),\theta_i(k))$ at $t$-th inner step of the SA scheme \eqref{SA2} with the initial  value $z_{i,1}(k) =x_i(k)$. Set  $x_i(k+1)   = z_{i, j_i(k) }(k).$}   Define $\xi_i(k)\triangleq ( \xi_{i,1}(k),\cdots, \xi_{i,j_i(k)}(k)),$ and 
  $\eta_i(k)\triangleq ( \eta_{i,1}(k) ,\cdots, \eta_{i,N_i(k)}(k) ),$  and   {$\mathcal{F}_k'\triangleq\sigma\{ x(0),i_l,\xi_{i_l}(l), \eta_{i_l}(l), d_{i_l,j}(l), 0 \leq l   \leq k-1, j\in \mathcal{N}\}$.}
  Then by Algorithm \ref{algo-randomized2},   $x(k)$ and $\theta(k)$  are adapted to $\mathcal{F}_k'$,
  and hence $y^i(k)$  is adapted to $\mathcal{F}_k$.  Thus,
     $T_i (y^i(k) ,\theta_i(k))$ is  adapted to $\mathcal{F}_k$  by its definition \eqref{def-br2}.
Then by Assumptions \ref{assp-miss} and  \ref{assp-diff},  we  obtain  the  same  bound  as
that of Lemma \ref{potential-lem-rate-sa}.
 Consequently,  Assumption \ref{assp-noise2}(b) is satisfied by
     setting $j_i(k)=\left \lfloor \Gamma_i(k)^{2(1+\delta)}\right \rfloor$, where $\delta>0$ and
 $\Gamma_i(k)$ is defined in Remark \ref{rem1}.

 	\subsection{ Convergence Analysis}
 We begin by proving a supporting Lemma,   an extension of the analogous determinstic (error-free) result from    \cite{bubeck2015convex}, which  will be used in the convergence  analysis of Algorithm \ref{algo-randomized2}.

\begin{lemma}\label{proj-lemma} Suppose Assumption    \ref{assp-miss} (c) holds. Let $\theta, y \in \Theta$ and suppose $\theta^+$ and $c_{\Theta}(\theta)$ are defined by . 
\begin{align} 
\theta^+ & := \Pi_{\Theta} \left( \theta - {1\over L_g} (\nabla_{\theta} g(\theta) + u)\right)  \mbox{ and }
 c_{\Theta}(\theta)  \triangleq L_g(\theta-\theta^+), \end{align} respectively.
Then the following holds.
\begin{align*}
	g(\theta^+)-g(y) & \leq c_{\Theta}(\theta)^T(\theta-y) - {1\over 2L_g} \|c_{\Theta}(\theta)\|^2 - u^T(\theta^+-y) - {\mu_g \over 2} \|\theta-y\|^2.
\end{align*}
\end{lemma}
{\em Proof.}
We begin by recalling the projection inequality
$$ \left(\theta^+- \left(\theta - {1\over L_g} \left(\nabla_{\theta} g(\theta) +  u \right)\right)\right)^T (\theta^+- y) \leq 0. $$
Consequently, we have that 
\begin{align}\label{eq1}
		\nabla_{\theta} g(\theta) ^T(\theta^+-y) \leq c_{\Theta}(\theta)^T(\theta^+-y) -u^T(\theta^+-y).
\end{align}
Then by the $\mu_g$-strong convexity and $L_g$-smoothness of $g(\cdot)$,  we may now derive the following   bound.
\begin{align*}
	g(\theta^+)-g(y) & = g(\theta^+)-g(\theta)+g(\theta) - g(y) \\ 		& \leq \nabla_{\theta} g(\theta)^T(\theta^+-\theta) + {L_g \over 2} \|\theta^+-\theta\|^2 + \nabla_{\theta} g(\theta)^T(\theta-y) - {\mu_g \over 2} \|\theta-y\|^2 \\
	& = 			\nabla_{\theta} g(\theta)^T(\theta^+-y) +{1\over 2L_g} \|c_{\Theta}(\theta)\|^2  - {\mu_g \over 2} \|\theta-y\|^2\\
	& \overset{\eqref{eq1}}{\leq} c_{\Theta}(\theta)^T(\theta^+-y) -u^T(\theta^+-y) + 
{1\over 2L_g} \|c_{\Theta}(\theta)\|^2  - {\mu_g \over 2} \|\theta-y\|^2\\
	& = c_{\Theta}(\theta)^T(\theta-y) - {1\over 2L_g} \|c_{\Theta}(\theta)\|^2 - u^T(\theta^+-y) - {\mu_g \over 2} \|\theta-y\|^2.  \qquad \Box
\end{align*} 

 \begin{theorem}[{\bf almost sure convergence for inexact BR with learning}] \label{thm-rand-2}~
Let  $\{ x(k)\}$ and $\{\theta(k)\} $ be generated by Algorithm \ref{algo-randomized2}.
Suppose  Assumptions   \ref{assp-miss}, \ref{assp-diff} and  \ref{assp-noise2}   hold. We further assume that
  $\mu$ used in \eqref{def-br2} satisfies
 $\mu>{L_x\over 2}+\sqrt{3}L_x\tau$,     and for every $i\in \mathcal{N}$,   $N_i(k)=\left \lfloor  \Gamma_i(k)^{2(1+\delta)}\right \rfloor$ for some $\delta>0$  \red{with $\Gamma_i(k)$  defined by \eqref{rem1}}. Then  \red{the sequences $\{\theta_i(k)\}$ and $\{x_i(k)\}$ satisfy the following.} 

\noindent (a)  For any $i\in \mathcal{N},$ $ \sum_{k=1}^{\infty} \| \theta_i(k)-\theta^*\|^2<\infty ~ a.s.,$ and
$ \sum_{k=1}^{\infty} \| \theta_i(k)-\theta^*\| <\infty ~a.s.$.

\noindent (b)
For any $i\in \mathcal{N}$,  $  \sum\limits_{k=0}^\infty    \|  T_i(y^i(k) ;\theta^*)- x_i(k) \|^2<\infty \quad a.s.$.

\noindent  (c)
 For almost all $\omega \in \Omega$, every limit point of  $x(k,\omega)$   is a Nash equilibrium.

\noindent  (d)  There exists a connected subset $X_c^* \subset X^*$ such that $d(x(k), X_c^*)  \xlongrightarrow [k \rightarrow \infty]{} 0~a.s.$.
\end{theorem}
 {\bf Proof.}
(a) For  $i_k=i,$ by    $e_{i,p}(k) $  defined in Assumption \ref{assp-noise2}(d)  we can  rewrite  \eqref{VSSG} as follows:
\begin{align}  \theta_i(k+1) &=\Pi_{\Theta}\left[\theta_i(k)- \red{{1\over L_g} \left(  \nabla  g(\theta_i(k))+u_k \right)}\right] ,\label{VSSG2}
\end{align}
where $\red{u_k}\triangleq\frac{1}{N_i(k)}  \sum_{p=1}^{N_i(k)}e_{i,p}(k) .$
\red{From Lemma~\ref{proj-lemma}, by setting $\theta = \theta_i(k)$, $u=u_k$, and $y = \theta^*$, we have that $\theta^{+}=\theta_i(k+1)$, $ c_{\Theta}(\theta_i(k)) = L_g(\theta_i(k)- \theta_i(k+1))$,   implying the following inequality. 
\begin{align}\notag
-c_{\Theta}(\theta_k)^T(\theta_i(k)-\theta^*)  & \leq -\underbrace{(g(\theta_i(k+1))-g(\theta^*))}_{\ \geq \ 0} - {1\over 2L_g}\|c_{\Theta}(\theta_i(k))\|^2 - {\mu_g \over 2} \|\theta_i(k)-\theta^*\|^2 -u_k^T (\theta_i(k+1)-\theta^*)\\
	& \leq  - {1\over 2L_g}\|c_{\Theta}(\theta_i(k))\|^2 - {\mu_g \over 2} \|\theta_i(k)-\theta^*\|^2 -u_k^T (\theta_i(k+1)-\theta^*). \label{bd-ct} 
\end{align}  
Since  $ c_{\Theta}(\theta_i(k)) = L_g(\theta_i(k)- \theta_i(k+1))$, we  may bound $\|\theta_i(k+1)-\theta^*\|^2$ as follows.  
\begin{align*}
&\|\theta_i(k+1)-\theta^*\|^2  = \|\theta_k - {1\over L_{g}} c_{\Theta}(\theta_i(k))   -\theta^*\|^2  = \|\theta_i(k)- \theta^*\|^2 +  {1\over L^2_g} \|c_{\Theta}(\theta_i(k))\|^2  -{2 \over L_g} c_{\Theta}(\theta_i(k))^T(\theta_i(k)-\theta^*)  \\&
\overset{\eqref{bd-ct}}{\leq}  \|\theta_i(k)- \theta^*\|^2 +  {1\over L^2_g} \|c_{\Theta}(\theta_i(k))\|^2  - {2 \over L_g}  \left(   {1\over 2L_g}\|c_{\Theta}(\theta_i(k))\|^2 + {\mu_g \over 2} \|\theta_i(k)-\theta^*\|^2 +u_k^T (\theta_i(k+1)-\theta^*)\right) \\
	& = \left(1-{\mu_g \over L_g}\right) \|\theta_i(k) - \theta^*\|^2   - {2 \over L_g} u_k^T (\theta_i(k+1)-\theta^*)\\
	& = \left(1-{\mu_g \over L_g}\right) \|\theta_i(k) - \theta^*\|^2   - {2 \over L_g} u_k^T (\theta_i(k+1)-\bar{\theta}_i(k+1) ) - {2 \over L_g} u_k^T (\bar{\theta}_i(k+1)  -\theta^*),
\end{align*}
where $\bar{\theta}_i(k+1) \triangleq  \Pi_{\Theta} \left[ \theta_i(k)- {1\over L_g}  \nabla_{\theta} g(\theta_i(k))  \right]$.
By \eqref{VSSG2}  and the non-expansivity of the Euclidean projector, one may obtain   that $ 
-u_k^T ( \theta_i(k+1) -\bar{\theta}_i(k+1) )    \leq \| u_k^T\|  \left \|  \theta_i(k+1) -\bar{\theta}_i(k+1) \right \| 
\leq \| u_k \|^2/L_g .$  Therefore, 
\begin{align*}
\|\theta_i(k+1) -\theta^*\|^2 \leq  \left(1-{\mu_g \over L_g}\right) \|\theta_i(k )  - \theta^*\|^2   +{2 \over L_g^2}\| u_k \|^2- {2 \over L_g} u_k^T (\bar{\theta}_i(k+1)   -\theta^*).
\end{align*}  
Taking expectations conditioned on $\mathcal{F}_k$ on both sides of the above equation,  we obtain the next inequality since $\theta_i(k)$ and $\bar{\theta}_i(k+1)$ are adapted to $\mathcal{F}_k$. 
\begin{align*}
\mathbb{E}[\|\theta_i(k+1) -\theta^*\|^2 \mid \mathcal{F}_k] & \leq \left(1-{\mu_g \over L_g}\right)
\|\theta_i(k)  - \theta^*\|^2   + {2\over L^2_g} \mathbb{E}[ \|u_k\|^2 \mid \mathcal{F}_k]     \quad {\scriptstyle  (\textrm{ by}~\mathbb{E}[  u_k  \mid \mathcal{F}_k] =0) }   
\\& \leq \left(1-{\mu_g \over L_g}\right) \mathbb{E}[\|\theta_i(k)-\theta^*\|^2] + \frac{2M_2^2}{L_g^2 N_i(k)},   \quad {\scriptstyle  (\textrm{ by~ Assumption~} \ref{assp-noise2}(d))}.
\end{align*}
} 
While for $i_k\neq i,$  $\mathbb{E}\left[ \|  \theta_i(k+1)-\theta^*\| \big | \mathcal{F}_k\right] = \| \theta_i(k)-\theta^*\|^2$.
\red{Then  by $\mathbb{P}(i_k=i)=p_i$, we obtain   that
\begin{align}\label{sub0}
\mathbb{E}\left[ \|  \theta_i(k+1)-\theta^*\|^2 \big | \mathcal{F}_k\right] \leq
 (1-   p_i  \mu_g / L_g  )  \| \theta_i(k)-\theta^*\|^2 +{2p_i  M_2^2 \over L_g^2 N_i(k) }.
\end{align} 
By   employing  the conditional  variant of  Jensen's inequality to \eqref{sub0},   we may conclude the following. 
\begin{align*}
\mathbb{E}\left[ \|  \theta_i(k+1)-\theta^*\|  \big | \mathcal{F}_k\right] &  \leq \sqrt{1-  p_i  \mu_g / L_g} \| \theta_i(k)-\theta^*\|+{\sqrt{2p_i}   M_2 \over  L_g \sqrt{N_i(k)} }
\quad{\scriptstyle(\textrm{by~} \sqrt{a^2+b^2}\leq  a+b \textrm{~for~} a,b\geq 0)}
\\& = \| \theta_i(k)-\theta^*\| -\left(1-\sqrt{1-   p_i  \mu_g / L_g}\right) \| \theta_i(k)-\theta^*\|+ {\sqrt{2p_i}  M_2 \over L_g \sqrt{N_i(k)} }.
\end{align*}} 
\red{Thus, by   using Theorem 1 in  \cite{robbins1985convergence},  and   $\sum_{k=1}^{\infty} 1/\sqrt{N_i(k)}<\infty, a.s. $ by  Remark \ref{rem1},}  we obtain that $ \sum_{k=1}^{\infty} \| \theta_i(k)-\theta^*\|<\infty ~ a.s.$, \red{and hence 
$ \sum_{k=1}^{\infty} \| \theta_i(k)-\theta^*\|^2<\infty ~ a.s.$.}

(b) Note that  $\nabla_{x_i} f_i(x;\theta^*)$    is Lipschitz continuous in $x\in X$ with Lipschitz constant $L_x$ by Assumption \ref{assp-miss}(b). Then similar to \eqref{potential-inequ2}  we obtain the following inequality for any $C>0$:
 \begin{equation} \label{potential-inequ2m}
\begin{split}  f_i\left(T_i(y^i(k),\theta^*),x_{-i}(k);\theta^*\right)+V_{k+1} &\leq f_i(x(k);\theta^*) +V_k
+ {L_x^2\tau^2  \over 2C } \| x(k+1)-x(k)\|^2  \\&-\left(\mu-{L_x+C \over 2}\right)\|T_{i}(y^i(k),\theta^*)-x_i(k)\|^2,
  \end{split}
\end{equation}
where $V_k\triangleq {L_x^2\tau  \over 2C }      \sum_{h=k-\tau+1}^k (h-k+\tau) \| x(h)-x(h-1)\|^2.$
  By  Assumptions   \ref{assp-diff}(a),   \ref{assp-diff}(b),    and  Jensen's inequality,  the following holds for any $x\in X$:
 \begin{align*}
 \|   \nabla_{x_i} f_i(x_i,x_{-i};\theta)\| & =\|\mathbb{E}[\nabla_{x_i} \psi_i(x_i, x_{-i};\theta;\xi )]\|
    \leq  \sqrt{\mathbb{E}[\|\nabla_{x_i} \psi_i(x_i, x_{-i};\theta;\xi )\|^2]} \leq M_1.
       \end{align*}
  Then     by the mean-value theorem   and Cauchy-Schwarz inequality,  we have that
  \begin{equation} \label{fb-term1m}
\begin{split}    & f_i\left(x_i(k+1),x_{-i}(k); \theta^*\right)- f_i\left(T_i(y^i(k),\theta^*),x_{-i}(k);\theta^*\right)
\\&= \left(x_i(k+1)-T_i(y^i(k),\theta^*)\right)^T\nabla_{x_i}   f_i\left(z_i(k+1),x_{-i}(k);\theta^*\right)
  \leq M_1\| x_i(k+1)-T_i(y^i(k),\theta^*)\|.
  \end{split}
\end{equation}      where $z_i(k+1)=\vartheta_{i,k}x_i(k+1)+(1-\vartheta_{i,k}) T_i(y^i(k),\theta^*)$ for some $\vartheta_{i,k}\in (0,1).$
By the triangle inequality, Lemma \ref{Lip-Tmap}, and \eqref{randomized2},    we may conclude that
\begin{align}\label{triag}
 \|x_i(k+1) -T_i(y^i(k),\theta^*)\|& \leq \|   T_i(y^i(k),\theta_i(k))-T_i(y^i(k),\theta^*)\|+ \|\varepsilon_i(k+1)\| \notag
\\& \leq L_t \|\theta_i(k)-\theta^*\|+ \|\varepsilon_i(k+1)\|.
\end{align}
\us{We then substitute  \eqref{triag} in \eqref{fb-term1m}  to}  obtain   the following bound:
\begin{equation} \label{bd-f1}
\begin{split}  &   f_i\left(x_i(k+1),x_{-i}(k); \theta^*\right)- f_i\left(T_i(y^i(k),\theta^*),x_{-i}(k);\theta^*\right)  \leq M_1L_t \|\theta_i(k)-\theta^*\|
+ M_1\|\varepsilon_i(k+1)\| .
  \end{split}
  \end{equation}
  Then     by  Algorithm  \ref{algo-randomized2},  Assumption \ref{assp-miss}(d), \eqref{potential-inequ2m} and \eqref{bd-f1},    we may obtain the following bound:
  \begin{align} \label{potent-inequ3}
 \notag \  &\quad P(x(k+1)  )- P( x(k)  )    = f_{i_k}( x_{i_k}(k+1),x_{-i_k}(k);\theta^*) -  f_{i_k}( x_{i_k}(k),x_{-i_k}(k);\theta^* )
  \\&  =  f_{i_k}\left(x_{i_k}(k+1),x_{-i_k}(k); \theta^*\right)- f_{i_k}\left(T_{i_k}(y^{i_k}(k),\theta^*),x_{-i_k}(k);\theta^*\right)
 \notag  \\& +f_{i_k}\left(T_i(y^{i_k}(k),\theta^*),x_{-i_k}(k);\theta^*\right)- f_{i_k}(x(k);\theta^*)   \notag  \\&  \leq
  M_1L_t \|\theta_{i_k}(k)-\theta^*\|+ M_1\|\varepsilon_{i_k}(k+1)\|+V_k-V_{k+1} \\ \notag&-\left(\mu-{L_x+C \over 2}\right)\|T_{i}(y^i(k),\theta^*)-x_i(k)\|^2+ {L_x^2\tau^2  \over 2C } \| x(k+1)-x(k)\|^2.   \end{align}
      By \us{the update rule in } Algorithm  \ref{algo-randomized2} and \eqref{triag},  we have that
      \begin{align*}
\| x(k+1)-x(k)\|&=\|x_{i_k}(k+1)-x_{i_k}(k) \|=
      \|  x_{i_k}(k+1)-T_{i_k}(y^{i_k}(k),\theta^*)+  T_{i_k}(y^{i_k}(k),\theta^*)-x_{i_k}(k)\|
    \\&\leq  L_t \|\theta_{i_k}(k)-\theta^*\|+ \|\varepsilon_{i_k}(k+1)\|+\|T_{i_k}(y^{i_k}(k),\theta^*)-x_{i_k}(k)\| ,
\end{align*}
\us{which when combined} with \eqref{potent-inequ3}  and $(a+b+c)^2\leq 3(a^2+b^2+c^2)$ \us{yields}  the following inequality:
   \begin{equation}\label{mis-inequ1}   \begin{split}   P(x(k+1) )+V_{k+1} &\leq P(x(k) ) +V_k+
M_1L_t \|\theta_{i_k}(k)-\theta^*\|+ M_1\|\varepsilon_{i_k}(k+1)\|  + {3L_x^2\tau^2  \over 2C } \|\varepsilon_{i_k}(k+1) \|^2 \\&-\underbrace{\left(\mu-{L_x+C \over 2}-{3L_x^2\tau^2  \over 2C }\right)}_{\triangleq \widetilde{C}}\|T_{i_k}(y^{i_k}(k),\theta^*)-x_{i_k}(k)\|^2+ {3L_x^2L_t^2\tau^2  \over 2C }
 \|\theta_{i_k}(k)-\theta^* \|^2 \\&\leq P(x(k) ) +V_k-\widetilde{C}\|T_{i_k}(y^{i_k}(k),\theta^*)-x_{i_k}(k)\|^2+
M_1L_t \|\theta_{i_k}(k)-\theta^*\| \\&   + {3L_x^2L_t^2\tau^2  \over 2C }
 \|\theta_{i_k}(k)-\theta^* \|^2  + M_1 \sum_{i=1}^N  \|\varepsilon_i(k+1)\|  + {3L_x^2\tau^2  \over 2C } \sum_{i=1}^N \|\varepsilon_i(k+1) \|^2  .
  \end{split}
\end{equation}
 Note that   $ x(k), V_k $,  $\theta_i(k),T_i (y^i(k);\theta^*)  ~\forall i\in \mathcal{N}$  are adapted to $\mathcal{F}_k$, and $i_k$ is independent of $\mathcal{F}_k$. Then by  taking expectations   conditioned on  $\mathcal{F}_k $ of \eqref{mis-inequ1}, and using   Corollary 7.1.2  in \cite{chow2012probability} and $\mathbb{P}(i_k=i)=p_i$,  it follows that 
%
   \begin{equation}     \label{mis-delay-cond-exp}
\begin{split}   \mathbb{E} \left[ P(x(k+1))+V_{k+1}  \big | \mathcal{F}_k \right] &  \leq   P(x(k)) +V_k - \widetilde{C} \sum_{i=1}^N p_i \|T_{i}(y^i(k),\theta^*)-x_i(k)\|^2
  \\&+M_1L_t  \sum_{i=1}^N p_i  \|\theta_i(k)-\theta^*\|       + {3L_x^2L_t^2\tau^2  \over 2C } \sum_{i=1}^N p_i  \|\theta_i(k)-\theta^* \|^2\\&+M_1 \sum_{i=1}^N  \mathbb{E} \left[    \| \varepsilon_i(k+1)\|        \big | \mathcal{F}_k\right]  + {3L_x^2\tau^2  \over  2C } \sum_{i=1}^N
    \mathbb{E} \left[    \| \varepsilon_i(k+1)\|^2        \big | \mathcal{F}_k\right]   .
  \end{split}
\end{equation}
By setting $C={\sqrt{3}L_x\tau\over 2}$ we derive ${C \over 2}+ {3L_x^2\tau^2  \over 2C }=\sqrt{3}L_x\tau .$
Thus, by taking  $\mu>{L_x\over 2}+\sqrt{3}L\tau$ it follows that  $ \widetilde{C} >0.$
Therefore,   by \red{ using Theorem 1 of \cite{robbins1985convergence}}, Assumption \ref{assp-noise2}(b) and  result (a), we have that $ \sum_{k=1}^{\infty} \sum_{i=1}^N p_i \|T_{i}(y^i(k),\theta^*)-x_i(k)\|^2<\infty ~a.s.$.
Then by $p_i\in (0,1)$ we obtain (b).

(c) The proof is \us{similar to} that of Theorem \ref{thm-rand-0}(b).

(d)  Since  $  \mathbb{E} \left[    \| \varepsilon_i(k+1)\|        \big | \mathcal{F}_k\right], \| \theta_i(k)-\theta^*\|, $
and $\|  T_i(y^i(k) ;\theta^*)- x_i(k) \|^2$ are  nonnegative    for  $k\geq 1, $
 by Assumption \ref{assp-noise2}(b), results (a) and (b)   we have the following for any $i\in \mathcal{N}$:
  \begin{equation*}
\begin{split} &  \sum_{k=1}^{\infty} \mathbb{E}[ \| \varepsilon_i(k+1)\|  ]=\mathbb{E} \left[\sum_{k=1}^{\infty}  \mathbb{E} \left[    \| \varepsilon_i(k+1)\|        \big | \mathcal{F}_k\right]  \right]< \infty,   \\&\sum_{k=1}^{\infty} \mathbb{E}[  \| \theta_i(k)-\theta^*\|  ]<\infty,
  \textrm{~and~} \sum_{k=1}^{\infty} \mathbb{E} \left[\|  T_i(y^i(k) ;\theta^*)- x_i(k) \|^2 \right]< \infty.
  \end{split}
  \end{equation*}
   Thus, by the Jensen's inequality we have  the following  for any $i\in \mathcal{N}$:
  \begin{equation}\label{limit_ti}
\begin{split} &  \lim_{k \to \infty} \mathbb{E}[ \| \varepsilon_i(k+1)\|  ]=0,~ \lim_{k \to \infty} \mathbb{E}[  \| \theta_i(k)-\theta^*\|  ]=0,
  \textrm{~and~}   \lim_{k \to \infty} \mathbb{E} \left[\|  T_i(y^i(k) ;\theta^*)- x_i(k) \|  \right]=0.
  \end{split}
  \end{equation}
Then by the  triangle inequality  and \eqref{triag}, we  obtain that
\begin{align*}
\mathbb{E}[\| x_i(k+1)-x_i(k)\|]& \leq  \mathbb{E}[ \|x_i(k+1) -T_i(y^i(k),\theta^*)\|]+ \mathbb{E} [\|x_i(k) -T_i(y^i(k),\theta^*)\|]
\\ & \leq  L_t  \mathbb{E}[\|\theta_i(k)-\theta^*\|]+\mathbb{E} [\|\varepsilon_i(k+1)\|]+ \mathbb{E} [\|x_i(k) -T_i(y^i(k),\theta^*)\|] \to 0, \us{ \mbox{ as $k \to \infty$.}}
\end{align*}
This  implies that $ \lim\limits_{k \rightarrow \infty } \mathbb{E}[\|x(k+1) -x(k)\|]  = 0.$
 The result follows  by proceeding as in  Theorem \ref{thm-rand-0}(c).
\hfill $\Box$

 \begin{remark} Note that in Algorithm  \ref{algo-randomized2},   we have utilized a variable sample-size stochastic approximation framework for updating each player's estimate of $\{\theta_i(k)\}$, \red{where each player  initiates its update from an independently selected  starting point.  The proposed scheme is similar to that developed in ~\cite{byrd2012sample} but we provide a distinct
proof, inspired by~\cite{bubeck2015convex}, through which the required
summability requirements are proven. Related schemes under differing assumptions with similar linear convergence rates have been studied over the last decade~(cf. \cite{shanbhag15budget,jofre19variance}). Accelerated variants of such schemes have also been studied in smooth~\cite{scholschmidt2011convergence} and nonsmooth~\cite{jalilzadeh2018optimal} regimes.  One might go a step further and note that players can utilize any update rule as long as it produces a sequence of iterates  satisfying  the requirement  that for every $i \in \mathcal N$, $ \sum_{k=1}^{\infty} \| \theta_i(k)-\theta^*\| <\infty ~a.s.$. } \end{remark}

Theorem \ref{thm-rand-2} shows that  the estimates of the equilibrium strategy
and the misspecified parameter generated by Algorithm  \ref{algo-randomized2}
converge almost surely to the set of Nash equilibria and to $\theta^*$, respectively.
Define the gap function  $G(x\us{;\theta^*}) \triangleq \sup_{y \in X} \nabla
P(x\us{;\theta^*})^T(x-y).$ The following result shows  the  convergence  in
mean of $x(k)$, characterized by the convergence of    $G(x(k))$ to zero  in the
mean sense.
\begin{theorem}[{\bf Convergence in mean of gap function.}]  \label{thm-rand-4}~
Let  $\{ x(k)\}$ and $\{\theta(k)\} $ be generated by Algorithm \ref{algo-randomized2}.
Suppose  Assumptions   \ref{assp-miss}, \ref{assp-diff} and  \ref{assp-noise2}   hold, and, in addition, that
    $\mu>{L_x\over 2}+\sqrt{3}L_x\tau$,      $\beta_i \in (0, 2 \mu_g/  L_g^2)$
and  $N_i(k)=\left \lceil  \Gamma_i(k)^{2(1+\delta)}\right \rceil$ for some $\delta>0$. Then
  $\lim\limits_{k \to \infty}\mathbb{E} \big[G(x(k)\us{;\theta^*})\big]= 0.$
\end{theorem}
{\bf Proof.} Similar to the proof of Theorem \ref{thm-rand-1},   we  have the following bound:    \begin{align*}
     & \mathbb{E}\left[  G(x(k)\us{;\theta^*}) \right]
 \leq  \sum_{i=1}^N\left(\mu D_{X_i} +M+L_x D_{X_i}\right) \mathbb{E}\left[\left \| T_i(y^i(k);\theta^*) -x_i(k)\right\|  \right]
 \implies 	\lim_{k \to \infty}   \mathbb{E}\left[  G(x(k)\us{;\theta^*}) \right]  \leq 0  \quad {\scriptstyle   \left(
				  \mathrm{by} ~ \eqref{limit_ti} \right) } .
   \end{align*}
   However, $G(x(k)\us{;\theta^*}) \geq 0$ since $x(k) \in X$, implying that $ \lim\limits_{k    \to \infty}     \mathbb{E}\left[  G(x(k)\us{;\theta^*}) \right]  = 0.$
\hfill $\Box$

\begin{remark}~ (i) Algorithm \ref{algo-randomized2} may   also  be  extended  to the generalized and weighted potential games with   misspecified parameters since Algorithm \ref{algo-randomized} is applicable to
the generalized and weighted potential games as shown in Section  \ref{Sec:IID}.
%

\noindent (ii)   The recent work by  \cite{jiang2017distributed} also considered the misspecified convex stochastic Nash games. This work was distinct in both its motivation and contributions.

\begin{itemize} \item  \cite{jiang2017distributed} consider monotone Nash games while we consider stochastic  potential games. Note that instances of one do not necesarily lie  within  the other.

\item  The update of equilibrium strategies in
 \cite{jiang2017distributed} utilizes projected gradient response while here we take inexact proximal best-response steps.

\item  The update for the misspecified parameter  in   \cite{jiang2017distributed}  is based on a  projected SG algorithm with a single  sampled gradient  per step and  with a decreasing step-size; here, we utilize an {\bf increasing sample-size}   projected SG scheme with a {\bf constant} step-size.
\end{itemize}

\end{remark}

\section{Preliminary numerics}
In this section, we \us{empirically} validate the performance of  Algorithm     \ref{algo-randomized}  and    Algorithm   \ref{algo-randomized2}
on the problem of congestion control and misspecified stochastic Nash-Cournot games, respectively.
\vspace{-0.1in}
\subsection{Congestion Control}
  We consider a congestion control problem on
a connected network characterized by a set of nodes $\mathcal{V} = \{1,\cdots, V\} $ and
a set of links $\mathcal{L} = \{1, \cdots, L\}$ connecting the nodes. There are
 $N$ users in the network,  where each  player $i$ aims at sending a   flow rate
 $x_i \in C_i=\{x_i \in \mathbb{R}: 0 \leq x_i \leq x_{i,\max}\}$
 from the source node $s_i$ to the destination node $d_i$  through a path $\mathcal{L}_i$ in the network.
 The upper bound $x_{i,\max}$ on   user $i$'s  flow rate might represent   a player-specific physical limitation.
 The  payoff function of player $i$    takes   as the difference of a player-specific
 pricing function and a utility function $U_i$ associated to the flow $x_i$  parameterized by
uncertainty $\xi_i,\zeta_i$:
 $$ \psi_i(x_i,x_{-i};\xi_i,\zeta_i) =\sum_{l \in \mathcal{L}_i} P_l\left(\sum_{j:l\in \mathcal{L}_j} x_j\right)-U_i\left(x_i, \xi_i,\zeta_i\right).$$
 The first term can be interpreted  as the price that  player  $i$ pays for the network resources with $P_l$ depending  on the
aggregated   flows on the link $l$.  Suppose that $P_l, l\in \mathcal{L}$ is convex and  $U_i ,i\in \mathcal{N}$ is  concave on $[0,x_{i,max}]$.
  Typical examples for the  pricing   and   utility functions  are   given by the following:
   \begin{align*}
 P_l=\frac{a_l}{b_l- \sum_{j:l\in \mathcal{L}_j} x_j }  \textrm{~and~}
 U_i=\xi_i \log(1+x_i+\zeta_i),
\end{align*}
where $\xi_i,\zeta_i$ are random variables.  Suppose each link $l\in \mathcal{L}$ in the network has a positive capacity $c_l$. Let us introduce a routing matrix $A\in \mathbb{R}^{L \times N}$,
where $[A]_{l,i}=1$ if $l \in \mathcal{L}_i$, and $[A]_{l,i}=0,$ otherwise.    The  capacity constraints of all links  can be expressed in the vector form as
  $Ax \leq c$ with $c=col\{c_l\}_{l=1}^L$.  For a fixed feasible $x_{-i}$,   we derive the  bound of  the  user    $i$'s flow rate $x_i$  denoted by
  $$0\leq X_i(x_{-i}) \leq \min_{l\in \mathcal{L}_i} \{c_l-\sum_{j\neq i} A_{l,j}x_j\}.$$
  The     $i$th user  aims at solving  the following  problem:
  \begin{align*}& \min_{x_i \in C_i \cap X_i(x_{-i})}\quad f_i(x_i,x_{-i})\triangleq \mathbb{E}\left[
\psi_i(x_i,x_{-i};\xi_i,\zeta_i) \right]. \end{align*}
Thus, the resulting problem is a generalized potential game  with the coupled constraint $X=\{x\in \mathbb{R}^n: Ax\leq c\}$ and the potential function  defined as follows:
  $$P(x)=\sum_{l \in \mathcal{L}} \frac{a_l}{b_l-\sum_{j:l\in \mathcal{L}_j} x_j}
  -\sum_{i\in \mathcal{N}} \mathbb{E}\left[U_i\left(x_i, \xi_i,\zeta_i\right)\right].$$
    Further, it is shown in Theorem 3.1 of \cite{alpcan2002game} that the congestion control  problem
    has a unique inner NE under appropriately chosen parameters.

\begin{figure}[!htb]
     \begin{minipage}{0.6\linewidth}
     \includegraphics[width=4.1in]{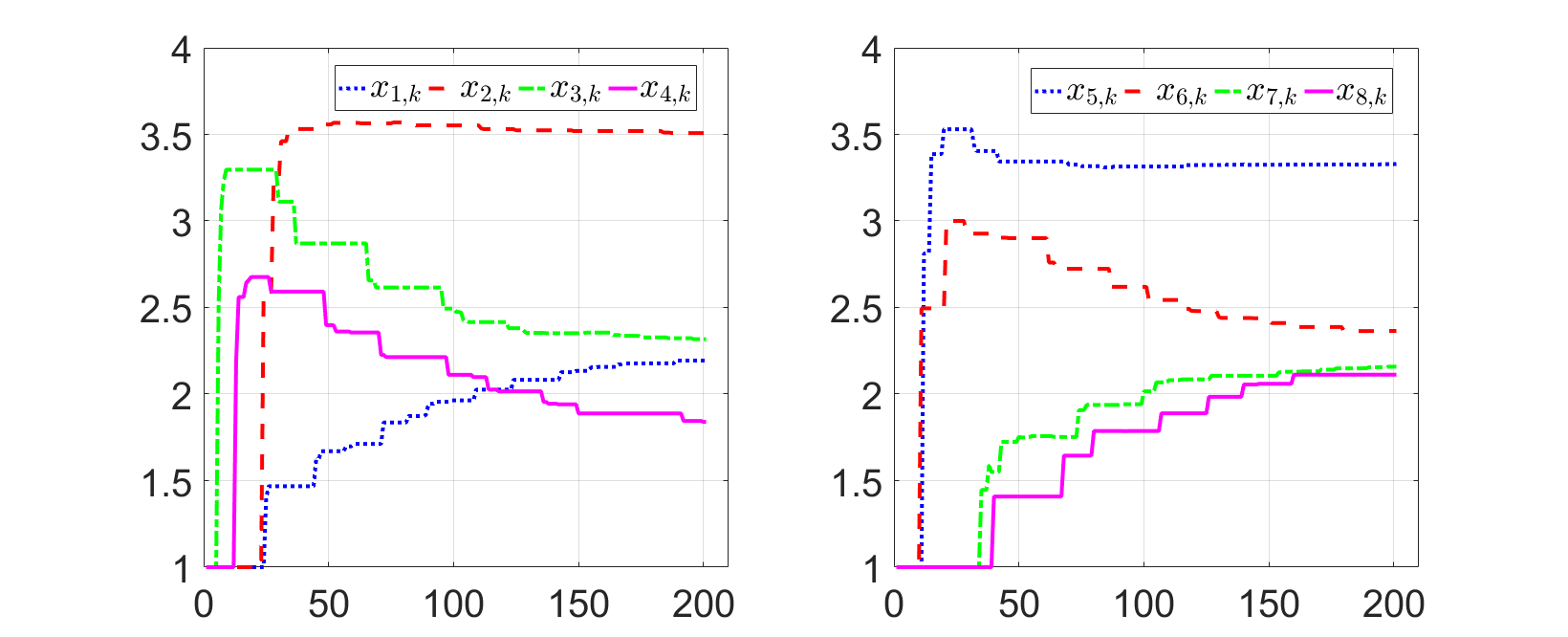}
  \caption{Flow rates of players (a single sample path)}    \label{fig-fr}
    \end{minipage}
   \begin{minipage}{0.4\linewidth}
       \centering
       \includegraphics[width=2.2in]{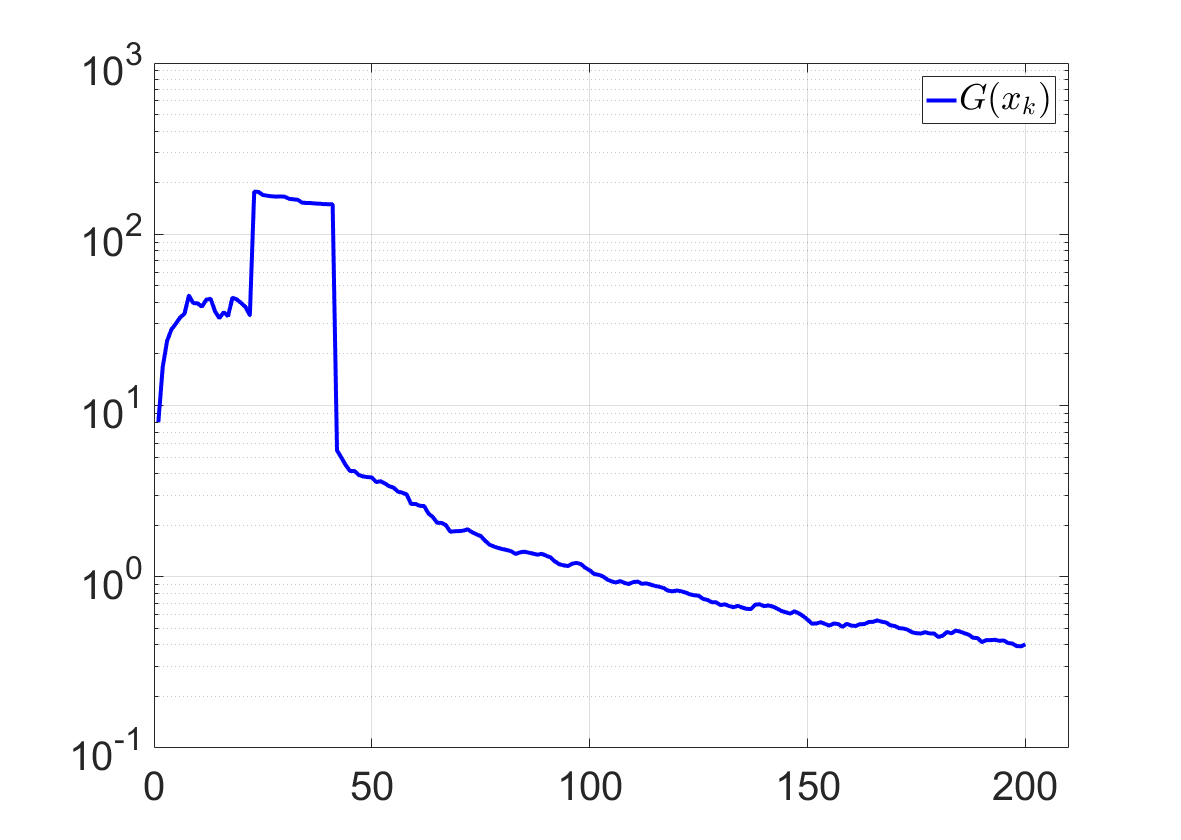}
      \caption{Trajectory of   $\mathbb{E}[G(x_k)]$}  \label{fig-opt}
    \end{minipage}%
 \end{figure}

 We \us{conducted}  numerical simulations for a network of $V=8$ nodes and $L=12$ links   shown in Figure \ref{Fig-One}.
 The parameters $a_l,b_l$   of  the utility function $P_l$ and the capacity constraint $c_l$  associated  to link $l \in \mathcal{L}$
 are given  in Figure \ref{Fig-One} as well. There are $N=8$ users  sending  flows  through  the network depicted    in  Figure \ref{Fig-One}.
The link paths of  user $i \in \mathcal{N}$ as well as   local parameters $\zeta_i,\xi_i,x_{i,\max},p_i$
 are given  in Table \ref{TAB-One}, where  $U[\tau_1,\tau_2]$  denotes the uniform distribution
  over the interval $[\tau_1,\tau_2]$. For any $k\geq0, i, j \in \mathcal{N}$,  the communication delays
$ d_{ij}(k)$ are independently generated from a uniform distribution on the  set  $\{0,1,\cdots, \tau\}$ with $\tau=4.$
  We carry out simulations for Algorithm   \ref{algo-randomized}, where the
  inexact solution \eqref{randomized-alg1} satisfying Assumption  \ref{assp-noise} are computed via the SA scheme (SA$_{i,k}$) with  $j_{i,k}= \lfloor   \Gamma_{i,k}^{3}\rfloor $ and   $\mu=1$.    The estimates of   each user‘s   equilibrium flow rates     are shown in   Figure  \ref{fig-fr},
  which demonstrates  the almost sure convergence of the iterates generated by Algorithm    \ref{algo-randomized}.
   Figure  \ref{fig-opt} displays  the     trajectory  of the mean gap function $\mathbb{E}[G(x_k)]$  calculated by averaging across 50 sample paths,
   which demonstrates      convergence in mean of the estimates generated by Algorithm    \ref{algo-randomized}.

   \begin{figure}[htbp]
     \centering
\begin{center}
\begin{tikzpicture}[->,>=stealth',shorten >=0.3pt,auto,node distance=1.8cm,thick,
  rect node/.style={rectangle, ball color={rgb:red,0;green,0.2;yellow,1},font=\sffamily,inner sep=1pt,outer sep=0pt,minimum size=14pt},
  wave/.style={decorate,decoration={snake,post length=0.1mm,amplitude=0.5mm,segment length=3mm},thick},
  main node/.style={shape=circle, ball color=green!20,text=black,inner sep=1pt,outer sep=0pt,minimum size=15pt},scale=0.75]
  \foreach \place/\i in {{(-2.1,3.1)/1},  {(-3,1.1)/2},
  {(-1.4,-1.4)/3},{(0.1,3.5)/4},  {(0.5,1)/5},
  {(1.8,-0.8)/6},{(3.5,2.9)/7}, {(3.8,0.1)/8}}
    \node[main node] (a\i) at \place {};

      \node at (-2.1,3.1){\rm \color{black}{$ 1$}};
      \node at (-3,1.1){\rm \color{black}{$2$}};
      \node at (-1.4,-1.4){\rm \color{black}{$3$}};
      \node at (0.1,3.5){\rm \color{black}{$ {4}$}};
      \node at (0.5,1){\rm \color{black}{${5}$}};
      \node at (1.8,-0.8){\rm \color{black}{${6}$}};
      \node at (3.5,2.9){\rm \color{black}{${7}$}};
        \node at (3.8,0.1){\rm \color{black}{${8}$}};

            \path[-,blue,thick]               (a1) edge (a4);
            \node at (-1.1, 3.6){\rm \color{black}{$ L_1$}};
                   \path[-,blue,thick]               (a4) edge (a7);
                          \node at (2.2, 3.4){\rm \color{black}{$ L_2$}};
                               \path[-,blue,thick]               (a7) edge (a8);

            \node at (4, 1.5){\rm \color{black}{$ L_3$}};

            \path[-,blue,thick]               (a6) edge (a8);
                \node at (2.8,-0.8){\rm \color{black}{$ L_4$}};
                        \path[-,blue,thick]               (a3) edge (a6);
                \node at (0.2,-1.3){\rm \color{black}{$ L_5$}};
            \path[-,blue,thick]               (a2) edge (a3);
 \node at (-2.5,0){\rm \color{black}{$ L_6$}};

            \path[-,blue,thick]               (a1) edge (a2);
             \node at (-3,2){\rm \color{black}{$ L_7$}};
                       \path[-,blue,thick]               (a1) edge (a3);
                        \node at (-1.4,0.4){\rm \color{black}{$ L8$}};

            \path[-,blue,thick]               (a5) edge (a3);

            \node at (4,1.5){\rm \color{black}{$ L_9$}};
            \path[-,blue,thick]               (a5) edge (a8);
                \node at (2.8,0.6){\rm \color{black}{$ L_{11}$}};

            \path[-,blue,thick]               (a5) edge (a7);
                              \node at (2,2){\rm \color{black}{$ L_{12}$}};
            \path[-,blue,thick]               (a1) edge (a5);
                \node at (-1,2){\rm \color{black}{$ L_{10}$}};
\end{tikzpicture}
\end{center}
\vskip 5mm
\scriptsize {
 \centering
     \begin{tabular}{|c|c|c|c|c|c|c|c|c|c|c|c|c|}
        \hline
 \diagbox[width=6.8em,trim=l]{Parameters}{Links $l$}  &1 &2 &3 &4 &5& 6 &7&8 &9 &10 &11 &12\\ \hline
     $a_l$& 5&4 &3& 5&4 &3& 5&4 &3& 5&4 &3\\ \hline
     $b_l$& 6&10 &8&6& 9& 5& 6& 5 &  6& 6& 8& 9\\ \hline
     $c_l$ & 5& 8 &6&  5& 8 & 4&  5& 4 & 4& 5& 7 & 8\\ \hline
      \end{tabular}
      }
   \caption{A network with 8 nodes and 12 links.}\label{Fig-One}

\vskip 5mm
\makeatletter\def\@captype{table}\makeatother
\newcommand{\tabincell}[2]{\begin{tabular}{@{}#1@{}}#2\end{tabular}}
 \scriptsize
 \centering
     \begin{tabular}{|c|c|c|c|c|c|}
        \hline
User $i$ & Link path & $\xi_i$ &$\zeta_i$  &$x_{i,\max}$ &$p_i$\\ \hline
     1 & $L_1,L_2,L_{12}$ & $U[10,12]$ & $U[0,1]$&3 &1/8\\\hline

  2 & $L_3,L_4,L_{5}$ &  $U[10,12]$ & $U[0,1]$&4&1/8\\ \hline

     3 & $L_{10},L_{11},L_{12}$ & $U[10,12]$& $U[0,1]$&4&1/8\\ \hline
   4 & $L_6,L_9,L_{12}$ &  $U[10,12]$ & $U[0,1]$&3&1/8\\ \hline

  5 & $L_5,L_8$ &   $U[10,12]$ & $U[0,1]$&5&1/8\\ \hline

     6 & $L_1,L_{2},L_7$ & $U[10,12]$ & $U[0,1]$&3&1/8\\ \hline

  7 & $L_3,L_{10},L_{11}$ &$U[10,12]$ & $U(0,1)$&4&1/8\\ \hline

8& $L_6 $ &  $U[10,12]$  & $U[0,1]$&3&1/8\\ \hline
      \end{tabular}
   \caption{Link paths and local parameters of all users}\label{TAB-One}
\end{figure}

 \subsection{Nash-Cournot Games  with Misspecified Parameters}  \label{sec-numerics-formulation}
We   apply   Algorithm   \ref{algo-randomized2}   to  the  networked Nash-Cournot game \cite{yu2017distributed,yi2017distributed}.
 ~ Suppose there are $N$ \us{firms}, regarded as the    set  of players  $\mathcal{N}=\{1,\dots,N\}$,
  competing over   $L$ markets denoted by $\mathcal{L}=\{1,\cdots,L\}$.   Firm  $i\in \mathcal{N}$ sells its products
    $x_i=(x_{i,1},\cdots, x_{i,n_i})\in \mathbb{R}^{n_i}$ to each  connected market with $n_i$ denoting the number of markets connected to firm $i$.
We use matrix $A_i \in  \mathbb{R}^{L \times n_i}$ to specify the participation of firm $i$ in the markets,
where     $[A_i]_{j,p}=1$   if firm  $i$ delivers its production $x_{i,p}$ to  market $j,$  and
$[A_i]_{j,p}=0,$ otherwise. The  production cost function of  firm $i$ is given by
$c_i(x_i;\xi_i)=(c_i+\xi_i)^T x_i $
for some positive parameter  $c_i\in \mathbb{R}^{n_i}$ and  random disturbance $\xi_i$ with mean zero.
Denote by $A=[A_1,\cdots,A_N],$ by $Ax=\sum_{i=1}^N A_ix_i\in \mathbb{R}^L$  and by  $S_j=[Ax]_j$ the  aggregated  products
of all connected firms  delivered to market $j$, where   $ [Ax]_j$ denotes the $j$-th entry of the vector $Ax.$
 Furthermore, the price of products sold in   market $j\in \mathcal{L}$ is assumed to follow a linear function corrupted by noise:
$$p_j(S_j;\zeta_j)=a_j^* +\zeta_j-  b_j^*S_j,$$
where $a_j^*>0,b_j^*>0$  are the pricing parameters, and the random disturbance $\zeta_j$ is zero-mean.
Then firm $i\in \mathcal{N}$ has a stochastic payoff function defined as follows:
 $$ \psi_i(x;\theta^*;\xi_i,\zeta_i) =c_i(x_i;\xi_i)-\sum_{j \in \mathcal{L} } p_j(S_j;\zeta_j) [A_ix_i]_j =
  (c_i+\xi_i)^T x_i- \Big(a^*+\zeta-B^* AX\Big)^T  A_ix_i ,$$
 where $a^*=col\{a^*_1,\cdots, a^*_L\},$  $\zeta=col\{\zeta_1,\cdots,\zeta_L\}$, $ B^*=diag\{b^*_1,\cdots, b^*_L\}, $
  and $\theta^*=(a^*,B^*)$  is unknown to all companies.
 Suppose    \us{firm} $i \in \mathcal{N}$  has \us{finite} production capacity $\us{X}_i=\{x_i\in \mathbb{R}^{n_i}:  0\leq x_i \leq \textrm{cap}_i\}.$  In this networked Cournot competition, firm $i  $ minimizes 
$c_i ^T x_i-  (A_ix_i)^T  a^*+(A_ix_i)^T   B^* \sum_{i=1}^N A_ix_i $ over $X_i$. If $P(x)$ is defined as $$P(x) \triangleq \sum_{i=1}^N c_i ^T x_i-  \left( \sum_{i=1}^N  A_ix_i \right)^T a^*
+ \Big(    col\left \{A_ix_i\right \}_{i=1}^N \Big)^T\chi \Big(    col\left \{A_ix_i\right \}_{i=1}^N \Big),$$
where $\chi={1 \over 2}(\mathbf{I}_N+\mathbf{J}_N) \otimes B^*$.
   Then  for any $ i\in \mathcal{N}$ {and for any $x_{-i} \in X_{-i}$,   equation  \eqref{def-potential} holds for all $x_i, x_i' \in X_i$. Thus,  the   Nash-Cournot game  admits  a  potential function  $P(x)$.
 By  definition of  $f_i(x_i,x_{-i};\theta^*)$, $\nabla_{x_i} f_i(x_i,x_{-i};\theta^*)$
  depends on $a_j^*, b_j^*,j\in \mathcal{L}$ if  \us{firm} $i$ sells  its products to  market $j $.
Each    firm $i $ can observe the historic data  about  the aggregated sales $S_j$  in market $j$
and  the price of products  $p_j=a_j^*+\zeta_j- b_j^*S_j$
 if the  market $j $ is   connected to firm $i$. As such, firm $i$ is able to learn the pricing parameters
  $a_j^*, b_j^*$ of the connected market $j$ through solving the following problem:
\begin{align}
\label{}
    &   \min_{a_j\geq 0,b_j\geq 0}  \mathbb{E} [( a_j-b_jS_j -p_j)^2]
\end{align}

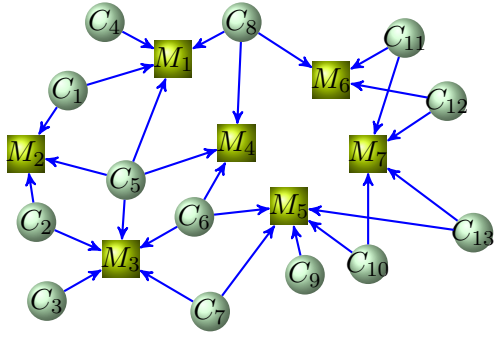
\begin{figure}
\begin{center}
\begin{tikzpicture}[->,>=stealth',shorten >=0.3pt,auto,node distance=1.8cm,thick,
  rect node/.style={rectangle, ball color={rgb:red,0;green,0.2;yellow,1},font=\sffamily,inner sep=1pt,outer sep=0pt,minimum size=14pt},
  wave/.style={decorate,decoration={snake,post length=0.1mm,amplitude=0.5mm,segment length=3mm},thick},
  main node/.style={shape=circle, ball color=green!20,text=black,inner sep=1pt,outer sep=0pt,minimum size=15pt},scale=0.7]


  \foreach \place/\i in {{(-3.2,2.2)/1},
  {(-3.8,-0.3)/2},
  {(-3.6,-1.8)/3},
  {(-2.5,3.5)/4},
  {(-2.1,0.5)/5},
  {(-0.8,-0.2)/6},
  {(-0.5,-2 )/7},
  {(0.1,3.5)/8},
  {(1.3,-1.3)/9},
  {(2.5,-1.1)/10},
  {(3.2,3.2)/11},
  {(4,2)/12},
  {(4.5,-0.5)/13}}
    \node[main node] (a\i) at \place {};

      \node at (-3.2,2.2){\rm \color{black}{$C_1$}};
      \node at (-3.8,-0.3){\rm \color{black}{$C_2$}};
      \node at (-3.6,-1.8){\rm \color{black}{$C_3$}};
      \node at (-2.5,3.5){\rm \color{black}{$C_4$}};
      \node at (-2.1,0.5){\rm \color{black}{$C_5$}};
      \node at (-0.8,-0.2){\rm \color{black}{$C_6$}};
      \node at (-0.5,-2 ){\rm \color{black}{$C_7$}};
      \node at (0.1,3.5){\rm \color{black}{$C_8$}};
      \node at (1.3,-1.3){\rm \color{black}{$C_9$}};
      \node at (2.5,-1.1){\rm \color{black}{$C_{10}$}};
       \node at (3.2,3.2){\rm \color{black}{$C_{11}$}};
      \node at(4,2) {\rm \color{black}{$C_{12}$}};
      \node at (4.5,-0.5){\rm \color{black}{$C_{13}$}};

  \foreach \place/\x in {{(-1.2,2.8)/1},{(-4,1)/2},{(-2.2,-1)/3},
    {(0,1.2)/4}, {(1,0)/5}, {(1.8,2.4)/6}, {(2.5,1)/7}}
  \node[rect node] (b\x) at \place {};

      \node at (-1.2,2.8){\rm \color{black}{$M_1$}};
      \node at (-4,1){\rm \color{black}{$M_2$}};
      \node at (-2.2,-1){\rm \color{black}{$M_3$}};
      \node at (0,1.2){\rm \color{black}{$M_4$}};
      \node at (1,0){\rm \color{black}{$M_5$}};
      \node at (1.8,2.4){\rm \color{black}{$M_6$}};
      \node at (2.5,1){\rm \color{black}{$M_7$}};

         \path[->,blue,thick]               (a1) edge (b1);
         \path[->,blue,thick]               (a1) edge (b2);

         \path[->,blue,thick]               (a2) edge (b2);
         \path[->,blue,thick]               (a2) edge (b3);

         \path[->,blue,thick]               (a3) edge (b3);

         \path[->,blue,thick]               (a4) edge (b1);

         \path[->,blue,thick]               (a5) edge (b1);
         \path[->,blue,thick]               (a5) edge (b2);
         \path[->,blue,thick]               (a5) edge (b3);
         \path[->,blue,thick]               (a5) edge (b4);

         \path[->,blue,thick]               (a6) edge (b3);
         \path[->,blue,thick]               (a6) edge (b4);
         \path[->,blue,thick]               (a6) edge (b5);

         \path[->,blue,thick]               (a7) edge (b3);
         \path[->,blue,thick]               (a7) edge (b5);

         \path[->,blue,thick]               (a8) edge (b1);
         \path[->,blue,thick]               (a8) edge (b4);
         \path[->,blue,thick]               (a8) edge (b6);

         \path[->,blue,thick]               (a9) edge (b5);

         \path[->,blue,thick]               (a10) edge (b5);
         \path[->,blue,thick]               (a10) edge (b7);

         \path[->,blue,thick]               (a11) edge (b6);
         \path[->,blue,thick]               (a11) edge (b7);

         \path[->,blue,thick]               (a12) edge (b6);
         \path[->,blue,thick]               (a12) edge (b7);

         \path[->,blue,thick]               (a13) edge (b5);
         \path[->,blue,thick]               (a13) edge (b7);
\end{tikzpicture}
\end{center}
    \caption{Networked Nash-Cournot: An edge from   $C_i$ to   $M_j$  implies  firm $C_i$  sells its products to market $M_j$.}\label{fig_network_cournot}
\end{figure}

\begin{figure}
  \centering
  \subfigure[~Learning  $a^*$]{
      \includegraphics[width=2in]{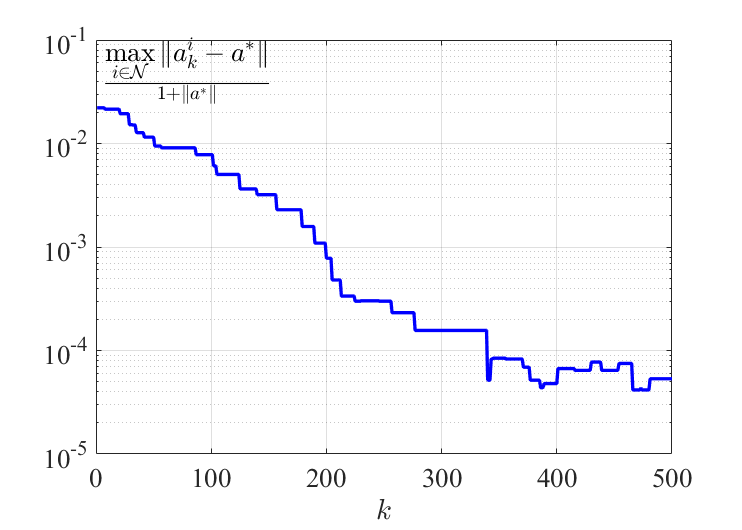} }
  \subfigure[~Learning $b^*$]{
    \includegraphics[width=2in]{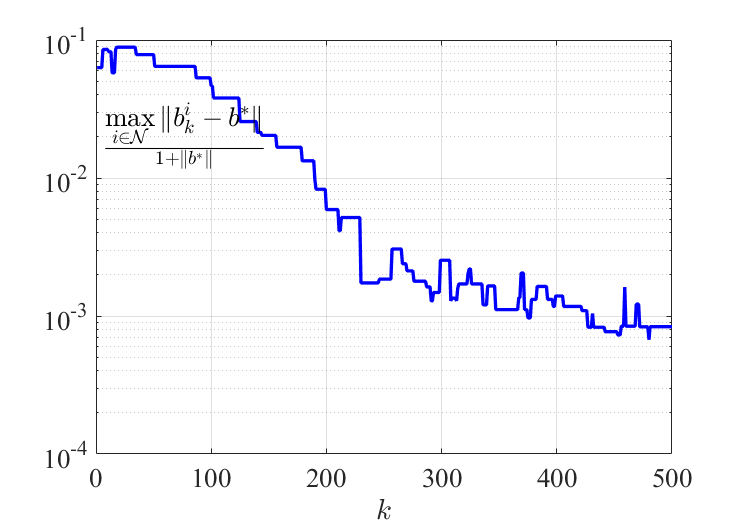} }
       \subfigure[~Learning $x^*$]{
    \includegraphics[width=2in]{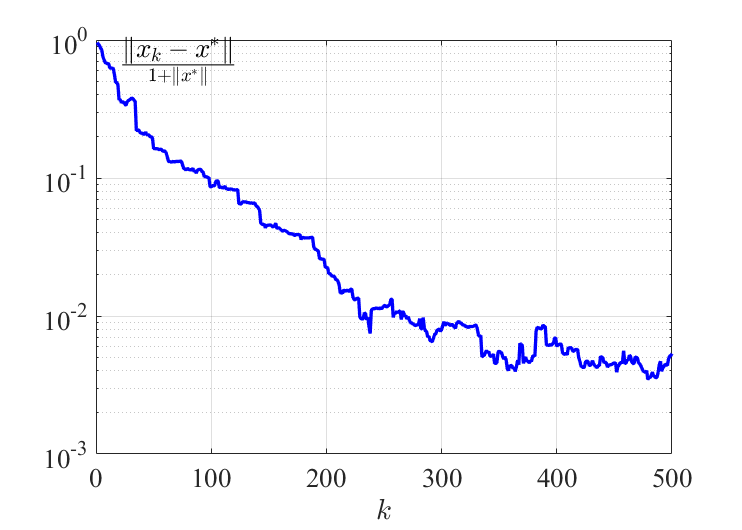}}
    \caption{The estimates of $a^*, b^*$ and $x^*$ (a single sample path)}   \label{fig-alg2}
\end{figure}

In the numerical investigation,   there are $V=13$ firms to sell their products to $L=7$ markets with the network shown in Figure \ref{fig_network_cournot}.
 Suppose  each component of $\textrm{cap}_{i } $ and the cost pricing parameter $c_i$
 of the  firm $i\in \mathcal{N}$  satisfy  uniform distributions specified by  $U[5,8]$ and $U[2,4]$.
The pricing parameters  $a_j^*,b_j^*$ of market $j\in \mathcal{L}$ are drawn from  uniform distributions $U[4,6]$ and $U[0.2,0.4]$, respectively.
  Suppose the random variables $\xi_i,i\in \mathcal{N}$ and $\zeta_j,j\in \mathcal{L}$
  are drawn from  uniform distributions $U[-c_i^*/8,c_i^*/8]$ and $U[-a_j^*/8,a_j^*/8]$, respectively.
  Suppose the historic aggregated sales  $S_j$ in  market $j\in \mathcal{L}$ satisfies  the uniform distribution  $U[0,5]$.
   For any $k\geq0, i, j \in \mathcal{N}$,  the communication delays
$ d_{ij}(k)$ are independently generated from a uniform distribution on the set  $\{0,1,\cdots, \tau\}$ with $\tau=4.$
  We carry out simulations for Algorithm   \ref{algo-randomized2}, where the
  inexact solution \eqref{randomized2} satisfying Assumption  \ref{assp-noise2} is computed via the SA scheme \eqref{SA2}
  with  $N_{i,k}=j_{i,k}= \lfloor   \Gamma_{i,k}^{3}\rfloor $,  $p_i=1/N~\forall i\in \mathcal{N},$ $\beta=0.1$ and   $\mu=5$.
  The scaled errors of learning schemes for the unknown parameters $a^*,b^*$ and the Nash equilibrium $x^*$ are provided in  Figs. \ref{fig-alg2},
  where $a_k^i$ and $b_k^i$ denotes the estimates of $a^*$ and $b^*$ given by firm $i$ at time $k.$
  The figure  demonstrates  the almost sure convergence of  Algorithm    \ref{algo-randomized2}.

   {\bf Comparison  with the asynchronous   SG  method:}  Suppose there are no communication delays among the players, i.e.,
  $\tau=0.$   Set $p_i=1/N~\forall i\in \mathcal{N},$ $\beta=0.1$ and  $\mu=5$.  Let
  $N_{i,k} = \lfloor   \Gamma_{i,k}^{3}\rfloor $ steps of the SA scheme \eqref{SA2}  be  taken at major iteration
  $k$ to obtain an   inexact solution to  \eqref{proximal},  where $\Gamma_{i,k}$ is defined in Lemma \ref{potential-lem-rate-sa}.
  Set $j_{i,k}= \lfloor   \Gamma_{i,k}^{3}\rfloor $  in equation \eqref{randomized2}.
   We  then carry out simulations for both  Algorithm \ref{algo-randomized2} and  the asynchronous SG algorithm,  which indeed is
    Algorithm \ref{algo-randomized2}  with equations \eqref{randomized2} and \eqref{VSSG}   replaced by
  \begin{align*}
x_{i,k+1} &=  \Pi_{X_i} \left[  x_{i,k}  -   \gamma_{i,k} \nabla_{x_i}   \psi_i( x_k,\theta_{i,k};\xi_{i,k})  \right] ,
\\ \theta_{i,k+1} &=\Pi_{\Theta}\left[\theta_{i,k }-\beta  \gamma_{i,k}  \nabla  g\left(\theta_{i,k},\eta_{i,k} \right)\right],
\end{align*}
where  $\gamma_{i,k}=\frac{1}{\Gamma_{i,k}^{0.6}}$.
We compare the two methods for the estimates of the equilibrium strategy $x^*$
 in terms of (i) the total number of the gradients steps (iteration complexity),    and (ii) the communication overhead  for achieving  the same accuracy.
Let   $K(\epsilon)$ denotes  the  smallest  total number of   SG steps  the players  have carried out to make
 $ \mathbb{E}\left[\frac{\|x_k-x^*\|}{1+\|x^*\|}\right]<\epsilon.$
 The empirically observed relationship between  $\epsilon$ and $K(\epsilon)$
for both methods are shown in  Figure \ref{fig-SG-BR}(a), where the  empirical
errors are calculated by averaging across 50 trajectories. From the figure,  it
is seen  that the  iteration complexity  are of the same orders while the
constant of SG is superior to that of Algorithm \ref{algo-randomized2}.  Since
SG requires the   rivals'  newest information to carry out a single  gradient
step, the  resulting  communication overhead   is proportional  to the total
number of gradient steps.  In contrast, Algorithm \ref{algo-randomized2}
carries  out multiple   gradient steps  without requiring the most recent
rivals'   information.  Further, the communication overhead  of   the two
methods are shown  in Figure \ref{fig-SG-BR}(b).  From the results
in Figure \ref{fig-SG-BR},   \us{upon termination},  Algorithm   \ref{algo-randomized2}
requires \us{approximately} about  10 times more gradient  steps than the  standard  SG
method while  \us{characterized by approximately} 500 times less  communication overhead.

\begin{figure}
  \centering
  \subfigure[~ Iteration Complexity]{
  \includegraphics[width=2.4in]{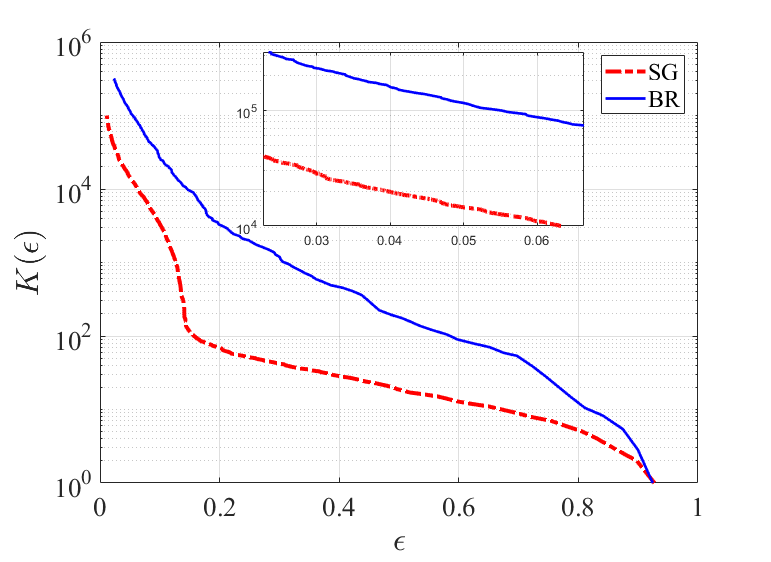} }
    \subfigure[~Communication Overhead]{
    \includegraphics[width=2.4in]{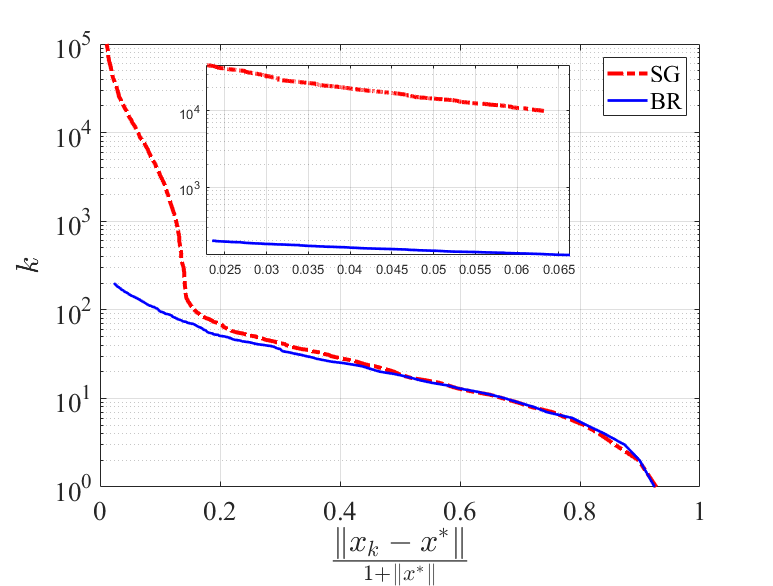}}
    \caption{Comparison of Stochastic Gradient Algorithm  and Inexact BR Algorithm  \ref{algo-randomized2}}   \label{fig-SG-BR}
\end{figure}

\section{Concluding Remarks}\label{sec:conclusion}
 This paper develops  an asynchronous inexact   proximal best-response   scheme
(\us{combined with joint learning})  to compute the Nash equilibrium of a
stochastic  potential  Nash game (\us{possibly corrupted by misspecification}).
\us{When player-specific problems are} convex, we show that the  estimates
generated by  the proposed schemes converge  almost surely  to a connected subset of
the NE set with uniformly bounded  delays. Since the game is
characterized by a possibly nonconvex potential function, the   schemes can
be viewed as randomized block coordinate descent schemes for a stochastic
nonconvex optimization problem which is block-wise convex. We   show
 that the gap function converges to zero  in mean for both schemes as well.
 {Furthermore, we prove almost sure convergence of the asynchronous  inexact BR scheme in the delay-free regime
 when the player-specific  problems are strongly convex.}
Finally,   we apply the developed methods to the congestion control problem and
the Nash-Cournot game, and demonstrate the  simulation results.

\def\cprime{$'$}

\end{document}